\documentclass[12pt]{amsart}

\usepackage[T1]{fontenc}
\usepackage{lmodern}
\usepackage{amsmath,amssymb,amsthm,mathtools, booktabs}
\usepackage{microtype}
\usepackage{needspace}
\usepackage[hidelinks]{hyperref}
\usepackage{xcolor}
\usepackage[normalem]{ulem}
\usepackage{soul}

\numberwithin{equation}{section}
\newtheorem{theorem}{Theorem}[section]
\newtheorem{lemma}[theorem]{Lemma}
\newtheorem{proposition}[theorem]{Proposition}
\newtheorem{corollary}[theorem]{Corollary}

\theoremstyle{definition}

\theoremstyle{remark}

\numberwithin{equation}{section}
\DeclareMathOperator{\reg}{reg}

\DeclareMathOperator{\Frac}{Frac}
\newcommand{\m}{\mathfrak m}
\newcommand{\Bn}{\overline B}
\newcommand{\rstabn}{\overline{\operatorname{rstab}}}

\begin{document}


\title{Regularity and stabilization of integral closures of powers of homogeneous ideals}

\author{Tran Nam Trung}
\address{Institute for artificial Intelligence, Vnu University of Engineering and Technology, Hanoi, Vietnam}
\email{tntrung@vnu.edu.vn}

\author{Guangjun Zhu}
\address{School of Mathematical Sciences, Soochow University, Suzhou, Jiangsu, 215006, P. R.~China}
\email{zhuguangjun@suda.edu.cn}


\begin{abstract}
In this paper, we establish explicit polynomial bounds for the regularity defect and the stabilization index of integral closures of powers of arbitrary homogeneous ideals over fields of any characteristic. Specifically, let $I \subseteq K[x_{1}, \dots, x_{n}]$ be a nonzero proper homogeneous ideal generated in degrees at most $d$. We prove that $\reg \overline{I^{s}} \le \rho(I)s + d^{(n-1)2^{n}}$ for every $s \ge 1$, where $\rho(I) \le d$ denotes the eventual slope of $\reg \overline{I^{s}}$. Furthermore, the stabilization bound is polynomial in $d$ for a fixed number of variables $n$. In characteristic zero, an upper bound for the first Hilbert coefficient yields a significantly sharper exponent. Our approach relies on exact equigenerated reductions, free lattice constructions within normalized Rees algebras, monic Frobenius presentations, and Tor-module regularity controls, complemented by generic source hyperplanes that reduce the stabilization behavior to the eventual nonvanishing of graded components.
\end{abstract}

\maketitle


\section{Introduction}
Cutkosky, Herzog and Trung \cite{CHT} and Kodiyalam \cite{Kod} proved that the regularity of powers of a homogeneous ideal is
eventually linear. The integral closures of powers also exhibit eventual linear regularity \cite{CHT,TW}. We give bounds for the defect
from the eventual slope and the onset of linearity in terms of the number of variables and the degrees of generators.

Let $R=K[x_1,\ldots,x_n]$ be a  standard graded  ring over a field $K$, and let  $ I\subset R$
be  a nonzero homogeneous ideal.  Let $d(I)$  be the maximum degree of all minimal homogeneous generators of $I$ and set
\[
\rho(I)=\min\{d(Q):Q\subseteq I
\text{ is a homogeneous reduction  of $I$}\}.
\]
Trung and Wang  \cite[Corollaries~3.3 and~3.4]{TW}  identified this number as the eventual slope for ordinary powers and integral
closures of powers. In the latter case,  the finiteness of the normalized Rees algebra  yields the eventual recurrence required in their result.
Write
\[
\reg\overline{I^s}=\rho(I)s+\bar b(I),  \text{ for $s\gg0$,}
\]
and define
\[
\begin{split}
	\Bn(I)&=\sup_{s\ge1}
	\{\reg\overline{I^s}-\rho(I)s\},\\
	\rstabn(I)&=\min\{s_0\ge1:
	\reg\overline{I^s}=\rho(I)s+\bar b(I)
	\text{ for all\ }s\ge s_0\}.
\end{split}
\]
We use the regularity conventions of \cite{Eis} and the integral closure terminology of \cite{HS}.

For a monomial ideal $I$, Hoa \cite[Theorem~3.13]{HoaIP}  proved the polynomial bound
\[
\rstabn(I)\le(n+1)(n+2)n^n d(I)^{2n^2}.
\]
Our contribution is to obtain polynomial bounds without the  monomial hypothesis. The proof uses an exact equigenerated
reduction and an independent stabilization argument in Section~\ref{sec:normalization}. For ordinary powers, Eisenbud and Ulrich    \cite{EU} provided stabilization results for an equigenerated ideal acting on a module $M$ equigenerated in a single degree,
under the assumption that $M/IM$ has a finite length. Bisui--H\`a--Thomas \cite{BHT} obtained bounds involving fiber invariants.

\begin{theorem}\label{thm:main}
Let $n \ge 3$ be an integer, let $0 \neq I \subsetneq K[x_1, \ldots, x_n]$ be a homogeneous ideal, and set $d=d(I)\ge2$.   Set
	\[
	M_n(d)=
	\begin{cases}
		d^{16},&n=3,\\[2mm]
		\displaystyle\left(\frac{(n-1)^2d^{2n}}8\right)^{2^{n-2}},&n\ge4.
	\end{cases}
	\]
	Then, over every field,  we have
	\[
	\max\{\Bn(I)+1,\rstabn(I)\}\le M_n(d).
	\]
	If $\operatorname{char}K=0$ and $n\ge4$, then
	\[
	\max\{\Bn(I)+1,\rstabn(I)\}
	\le\min\left\{M_n(d),
	\left(\frac{n(n-1)^2d^{2n-1}}4\right)^{2^{n-2}}\right\}.
	\]
\end{theorem}

For a fixed $n$, these bounds are polynomial in $d$ and independent of the number of generators. For $n\ge4$, the characteristic-zero estimate has degree $(2n-1)2^{n-2}$ in $d$,  compared to $2n2^{n-2}$ for the  arbitrary-field estimate.

Since integral closure can increase regularity \cite{Misra}, degree bounds for $\reg I$ do not directly bound $\reg\overline I$.
The next corollary bounds $\reg\overline I$ in terms of the generating degrees of $I$,  analogous to the bounds for $\reg I$ by  Caviglia and Sbarra  \cite{CS} and Hoa \cite{HoaReg}.
It applies uniformly to all powers, with the eventual slope $\rho(I)$ as the coefficient of $s$.
\begin{corollary}\label{cor:regularity-powers}
	Let $0\ne I\subset K[x_1,\ldots,x_n]$ be a homogeneous ideal whose minimal generators have degrees at most $d$.
	Then, for every $s\ge1$,
	\[
	\reg\overline{I^s}\le \rho(I)s+d^{(n-1)2^n}.
	\]
\end{corollary}

Our proof uses a bound for normalization in its inherited grading, which depends only on dimension and multiplicity.

\begin{theorem}\label{thm:normalization}
	Let $A$ be a geometrically integral standard graded algebra over  a field $K$, of dimension $r\ge3$ and multiplicity $e$.
	Let $C$ be its normalization in $\operatorname{Frac}A$,  endowed  with the inherited grading. Then
	\[
	\reg C\le e^{2^{r-1}}.
	\]
\end{theorem}
Here $\reg C$ is computed with respect to $A_+$, or equivalently with respect to a linear Noether normalization after a purely
transcendental field extension. For comparison, Hoa \cite[Theorem~1.5(i)]{HoaReg} established
\[
\reg S\le\left[\frac{e(e+1)}2\right]^{2^{r-2}}-1
\]
for a standard graded domain $S$ of dimension $r\ge2$ and multiplicity $e$.  Because the inherited grading on $C$ need not be
standard,  Hoa's bound does not apply directly to $C$.
Earlier work measures normalization complexity through generator degrees \cite{UV}, Hilbert coefficients and normalization chains
\cite{PV}, and normalization indices of ideals \cite{PUVV}. These indices concern  the generation of the normalized Rees algebra;
$\rstabn(I)$ records the onset of linearity of its strand regularities.

The rest of this paper is organized as follows. In Section  \ref{sec:equigenerat}, we replace $I$ with an equigenerated ideal $J=(I_{\rho})$, which preserves $\reg\overline{I^{s}}$ for all $s\ge1$. Section \ref{sec:normalization}  is devoted to the proof of the normalization bounds, free lattice constructions, Frobenius presentations, and the stabilization argument via generic source hyperplanes and graded components.

\section{Equigenerated reduction and degree}
\label{sec:equigenerat}

In this section, we reduce the study of $\reg\overline{I^s}$ to ideals generated in a single degree and record the saturation and elementary cases needed below.

Write $\m=(x_1,\ldots,x_n)$. For a homogeneous ideal $J$, let $(J_a)$ be the ideal generated by its degree-$a$ component. Unless stated otherwise, field extensions are purely transcendental. They preserve regularity and commute with normalization, since polynomial extension and localization preserve normality. Thus, we may choose generic linear forms over an infinite field.

\begin{lemma}\label{lem:truncation}
If $Q\ne0$ is  a homogeneous ideal  generated in degrees at most $a\ge1$, then  for all $s\ge1$,
\[
 \overline{(Q_a)^s}=(\overline{Q^s})_{\ge as},\qquad
 \reg\overline{(Q_a)^s}=\max\{\reg\overline{Q^s},as\}.
\]
\end{lemma}
\begin{proof}
Factoring homogeneous   components gives $(Q_a)^s=(Q^s)_{\ge as}$. If $f\in\overline{Q^s}$ is a  homogeneous  element of degree $b\ge as$, then  it satisfies an integral dependence equation with coefficients $c_i \in Q^{si}$.
By taking the homogeneous components of the equation, we can choose the coefficients $c_i$ to be homogeneous of degree $ib$. Since $Q$
is generated in degrees at most $a$, the homogeneous component of degree $ib$ of $Q^{si}$ coincides with the corresponding truncation of $(Q_a)^{si}$. Thus, $c_i \in (Q_a)^{si}$ because $ib \ge asi$. This implies that $f$ satisfies an integral equation over $(Q_a)^s$, and hence $f \in \overline{(Q_a)^s}$. This holds for all homogeneous components of degree $b \ge as$, giving $(\overline{Q^s})_{\ge as} \subseteq \overline{(Q_a)^s}$.
Conversely, since $(Q_a)\subseteq Q$ and $(Q_a)^s\subseteq\m^{as}$, and $\m^{as}$ is integrally closed, we have $\overline{(Q_a)^s}\subseteq\overline{Q^s}\cap\m^{as}=(\overline{Q^s})_{\ge as}$.

For the  equality concerning the regularity, for any non-zero ideal $J$ with $J_b \ne 0$, the standard truncation exact sequence gives $\reg J_{\ge b}=\max\{\reg J,b\}$. Applying this to $J=\overline{Q^s}$ and $b=as$ yields $\reg\overline{(Q_a)^s} = \reg(\overline{Q^s})_{\ge as} = \max\{\reg\overline{Q^s},as\}$.
\end{proof}

\begin{lemma}\label{lem:lower}
	If $Q$ is  a  nonzero proper homogeneous ideal generated in degrees at most $\rho=\rho(Q)$, then $\reg\overline{Q^s}\ge\rho s$ for every $s\ge1$.
\end{lemma}
	\begin{proof}
		Let $A = R[Qt]$ and $C = \bigoplus_{s \ge 0} \overline{Q^s}t^s$.  By lying over for $A\subseteq C$, every prime of $A$ containing $\mathfrak m A$ is the contraction of a prime of $C$ containing $\mathfrak m C$. Therefore the kernel of $A/\mathfrak m A\to C/\mathfrak m C$ consists of nilpotent elements and is nilpotent by Noetherianity.
		
	We claim that some generator $f \in Q_\rho$ has a nonnilpotent class in $A/\mathfrak{m}A$: otherwise, modulo $\mathfrak{m}Q^k$, every sufficiently long product of generators of $Q$ would contain a generator of degree less than $\rho$. By graded Nakayama's lemma, this would imply $Q^k = (Q_{\le \rho-1})Q^{k-1}$ for $k \gg 0$, which directly contradicts the definition of $\rho$.
		
		Consequently, the class of $ft$ stays nonnilpotent in $C/\mathfrak{m}C$. Therefore, $f^s \notin \mathfrak{m}\overline{Q^s}$ for every $s \ge 1$, and by Nakayama's lemma again, $\overline{Q^s}$ has a minimal generator of degree $\rho s$, which forces $\reg\overline{Q^s} \ge \rho s$.
	\end{proof}

\begin{proposition}\label{thm:reduction}
	Let $ I \subsetneq R$ be a nonzero  homogeneous ideal, and put $\rho = \rho(I)$ and $J = (I_\rho)$. Then
	\[
	\reg \overline{I^s} = \reg \overline{J^s} \quad \text{ for all $s\ge1$}.
	\]
	Consequently, $\overline{B}$ and $\overline{\mathrm{rstab}}$ agree for $I$ and $J$.
\end{proposition}
\begin{proof}
Since the ideal $Q = (I_{\le \rho})$ contains a reduction of $I$, it is itself a reduction of $I$. Consequently, $\overline{Q^s} = \overline{I^s}$ by \cite{HS}.  A reduction of $I$ generated in degrees at most $\rho$ is contained in $Q$ and also reduces $Q$, so $\rho(Q)\le\rho$. Conversely, every reduction of $Q$ reduces $I$ by transitivity, giving $\rho(Q)\ge\rho$. Given that $Q_\rho = I_\rho = J_\rho$, the two preceding lemmas together prove the claim.
\end{proof}

\begin{lemma}[High-degree saturation]\label{lem:normal-saturation}
	If $I$ is a nonzero homogeneous ideal generated in degrees at most $a$, then
	\[
	((\overline{I})^{\mathrm{sat}})_j = (\overline{I})_j \quad  \text{ for all $j\ge  a$}.
	\]
	In particular, if $L$ is generated in degree $\rho$, then $\overline{L^s} = ((\overline{L^s})^{\mathrm{sat}})_{\ge \rho s}$.
\end{lemma}
\begin{proof}
We employ the valuative argument outlined in \cite[Lemma 3.6]{Misra}. Let $f$ be a homogeneous element in $(\overline{I})^{\mathrm{sat}}$ of degree $j \ge a$, and let $v$ be a discrete valuation which is nonnegative on $R$. Choose a variable $x_h$ having the minimum value among all variables, and write $I = (g_1, \dots, g_q)$ with $a_i = \deg g_i \le a$
	
	Consider the  ring $K[x_i/x_h : i \neq h]$ and set $I' = (g_i/x_h^{a_i})$.  A power of $x_h$ times $f$ lies in $\overline I$, so $f/x_h^j$ is integral over the extension of $I'$ to $K[x_i/x_h:i\ne h][x_h,x_h^{-1}]$. Taking the degree-zero part of a monic integral equation gives $f/x_h^j\in\overline{I'}$. Thus, the valuation inequality yields
	\[
	v(f) - jv(x_h) \ge \min_i \{v(g_i) - a_i v(x_h)\}.
	\]
	At an index attaining the minimum, since $j \ge a \ge a_i$, it follows that $v(f) \ge v(g_i) \ge v(I)$. Therefore, the valuative criterion directly proves the first assertion. Finally, the second assertion is established by additionally noting that $\overline{L^s} \subseteq \mathfrak{m}^{\rho s}$.
\end{proof}

\begin{proposition}\label{prop:elementary-cases}
	 Let $I\subsetneq R$ be a nonzero homogeneous ideal. Suppose that at least one of the following holds: $n\le2$, $\rho(I)=1$, or $I$ is $\m$-primary. Then
	\[
	\reg\overline{I^s}=\rho(I)s\quad    \text{ for all     } s\ge 1,\qquad
	\Bn(I)=0,\qquad\rstabn(I)=1.
	\]
\end{proposition}
\begin{proof}
By Proposition~\ref{thm:reduction}, we can replace $I$ with the ideal $L=(I_\rho)$ generated by the degree-$\rho$ elements of $I$, where $\rho=\rho(I)$.
Since $L$ is a reduction of $I$, we have $\overline{I^s}=\overline{L^s}$ for every $s\ge1$, so it suffices to prove the assertion for $L$.
If $\rho(I)=1$, the ideal $L$ is generated by linear forms. Powers of ideals generated by linear forms are integrally closed and have linear free resolutions, which implies $\reg\overline{I^s} = s = \rho(I)s$.

If $I$ is $\mathfrak{m}$-primary, then so is $L$: indeed, since $L$ is a reduction of $I$, we have $\sqrt{L}=\sqrt{I}=\mathfrak m$. If $L$ is $\mathfrak{m}$-primary, then its generators satisfy $L\mathfrak{m}^{\rho k}=\mathfrak{m}^{\rho(k+1)}$ for $k\gg0$,
which implies that $L$ reduces $\mathfrak m^\rho$, and thus
\[
\overline{L^s}=\overline{(\mathfrak m^\rho)^s}=\mathfrak m^{\rho s}.
\]
Consequently, $\reg\overline{I^s}=\rho s$. In the case $n=1$,  $I$ is a principal ideal, and the assertions are trivial.
In the case $n=2$, we can write $L$ as $L = fJ$, where  $f$  is the greatest common divisor of the generators of $L$, and the generators of $J$ have no common factor of positive degree. Thus $f$ is homogeneous, and  $J = R$ or $J$ is $\mathfrak{m}$-primary. Since integral closure commutes with principal factors in polynomial rings, we have $\overline{L^s} = f^s \overline{J^s}$.
If $J=R$, then $\overline{L^s}=(f^s)$ and hence $\reg\overline{L^s}=s\deg f=\rho s$. If $J$ is $\mathfrak m$-primary, the preceding argument gives $\overline{J^s}=\mathfrak m^{ds}$, where $d=\deg J=\rho-\deg f$. Therefore
	\[
	\reg\overline{L^s}=\reg(f^s\mathfrak m^{ds})=s\deg f+ds=\rho s.
	\]

	Combining these cases, the linear growth of regularity $\reg\overline{I^s} = \rho s$ follows, along with $\Bn(I) = 0$ and $\rstabn(I) = 1$.
\end{proof}

Henceforth, assume that $I$ is 	an ideal generated in degree $\rho \ge 2$. Give $A = R[It]$ the bigrading
\[
A_{(u,s)} = (I^s)_{\rho s + u} t^s.
\]
Thus $\deg x_i = (1,0)$, $\deg(ft) = (0,1)$ for $f \in I_\rho$, and total degree is $u+s$. The normalization $C = \bigoplus_{s \ge 0} \overline{I^s} t^s$ inherits this nonnegative bigrading. It is geometrically integral and has dimension $n+1$.
Unless specified otherwise, regularity refers to the total grading. Write $\ell(I) = \operatorname{dim} K[I_\rho t]$ for the analytic spread.

\section{Free lattices and normalization}
\label{sec:normalization}

In this section, we develop the technical machinery regarding free lattices, Frobenius presentations, and cone decompositions needed to prove our main normalization and stabilization theorems.

All polynomial parameters have degree one. We use the regularity inequalities for short exact sequences \cite{Eis, CFN}. For an algebra finite over a linear Noether normalization $P$, regularity can be computed over $P$, since the positive ideals have the same radical. A free lattice is a graded free $P$-submodule of full rank, not necessarily an algebra.

\subsection{General normalization bounds}
\label{subsec:general-normalization}
To prove Theorem~\ref{thm:normalization}, we construct free lattices, compare cone degrees, and recover the normal algebra by Frobenius.

\paragraph{\textbf{Free lattices and cone comparison.}}
We first construct a graded free $P$-submodule of $F$ with explicitly bounded basis degrees, closely related to \cite[Proposition 5.3]{DV}, for subsequent cone degree comparisons.

\begin{lemma}[A bounded free lattice]
	\label{lem:bounded-field-lattice}
Let $P \subseteq F \subseteq V$ be graded domains finite over the standard graded polynomial ring $P$, with $\operatorname{Frac} F = \operatorname{Frac} V$. Suppose $F$ is  generated as a $P$-algebra in degrees at most $B$, and put $r_0 = \operatorname{rank}_P V$. Then $F$ contains a graded free $P$-module
\[
U = \bigoplus_{j=0}^{r_0-1} P(-d_j), \quad d_0 = 0, \quad 0 \le d_j \le jB.
\]
Thus $\max d_j \le B(r_0 - 1)$ and $\sum_j d_j \le Br_0(r_0 - 1)/2$.
\end{lemma}
\begin{proof}
Let $K_0 = \operatorname{Frac} P$. Since $V$ is a finite $P$-module, $V \otimes_P K_0$ is a field of dimension $r_0$ over $K_0$. Choose homogeneous generators of $F$ as a $P$-algebra. For $j \ge 0$, let $W_j$ be the $K_0$-vector space spanned by monomials of length at most $j$, giving a filtration $K_0 = W_0 \subseteq W_1 \subseteq W_2 \subseteq \cdots$.

If $W_j = W_{j+1}$ for some $j$, then multiplication by the generators preserves $W_j$. Since $1 \in W_j$, it contains all monomials in these generators, yielding $F \otimes_P K_0 \subseteq W_j$. As $\operatorname{Frac} F = \operatorname{Frac} V$, we have $F \otimes_P K_0 = V \otimes_P K_0$, so $W_j = V \otimes_P K_0$.

Since $\dim_{K_0}(V \otimes_P K_0) = r_0$ and dimensions strictly increase until stabilization, we can choose homogeneous monomials $\beta_j \in W_j$, for $j=0, \dots, r_0-1$ with $\beta_0 = 1$, that are linearly independent over $K_0$ (and hence over $P$). Setting $d_j = \operatorname{deg} \beta_j \le jB$, their $P$-span is the graded free module $U = \bigoplus_{j=0}^{r_0-1} P(-d_j)$. The bounds on $\max_j d_j$ and $\sum_j d_j$ follow immediately.
\end{proof}

\begin{lemma}[The first Hilbert coefficient]
\label{lem:normal-finite-e1}
Let $V$ be a normal graded $K$-domain generated in nonnegative degrees, with $V_0=K$ and dimension $d \ge 2$. If $V$ is finite over the polynomial subring $P=K[z_1,\ldots,z_d]$ with $\deg z_i=1$, then its first Hilbert coefficient is nonnegative, i.e., $e_1(V)\ge0$.
\end{lemma}
\begin{proof}
Let $A$ be the local ring $V_{V_+}$ and $J=(z_1,\ldots,z_d)A$ be the ideal generated by the parameters. For any $s \ge 1$, the integral closure of $J^s$ satisfies
 \[
\overline{J^s}=V_{\ge s}A.
\]
To see this, note that each homogeneous element $v \in V_t$ is integral over $J^t$ due to its integral dependence equation over $P$. Conversely, the valuation defined by the least homogeneous degree shows that integrality over $J^s$ forces elements to lie in $V_{\ge s}A$.
Consequently, the length of the quotient ring can be computed as
\[
\operatorname{length}_A(A/\overline{J^s}) = \sum_{t=0}^{s-1}\dim_KV_t,
\]
 which implies that the extended Hilbert coefficient matches $e_1(V)$, namely $\overline e_1(J)=e_1(V)$.

Finally, since $A$ is an excellent normal local ring, its completion remains normal, allowing us to conclude by applying \cite[Theorem~3.1]{MQ}.
\end{proof}

\paragraph{\textbf{Cone decompositions.}}
We recall the notation needed from \cite{Liang,Nagel}.  References to the result numbers in \cite{Liang} use its arXiv version. In a graded free module over $T=k[w_1,\ldots,w_r]$, a cone is a subspace of the form $h\,k[U]$, where $h\ne0$ is
homogeneous and $U\subseteq\{w_1,\ldots,w_r\}$. Its degree and dimension are defined to be $\deg h$ and $|U|$, respectively. A cone decomposition is a finite
direct sum of cones. It is called $q$-standard if every positive-dimensional cone has degree at least $q$, and whenever a cone of degree $a$ and dimension $i>0$ occurs, there is also a cone of dimension at least $i$ in each degree $q,\ldots,a$. It is called $q$-exact if, in addition, all positive-dimensional cones have distinct degrees.
For a $q$-exact decomposition $\mathcal Q$, define
\[
 b_i(\mathcal Q)=
 \max\bigl(\{q\}\cup
 \{1+\deg C:C\in\mathcal Q,\ \dim C\ge i\}\bigr).
\]
For a fixed $q$, these constants are uniquely determined by the Hilbert function of $\mathcal Q$
\cite[Lemma~26]{Liang}.

For a nonzero graded module $N$ generated in degree zero, we denote by $b_0(N)$ the constant associated with a $0$-exact decomposition of its
normal-form space, as in \cite[Proposition~2.10]{Nagel},
and $\beta_i(D)$ represents the constants with auxiliary parameter $D$, distinguishing them from those of $N$.

\begin{lemma}[Comparison with a free lattice]
\label{lem:free-lattice-exact}
Let $T=K[w_1,\ldots,w_r]$ with $r\ge3$, and let $\mathcal M\subseteq T^m$ be generated in degrees at most $D\ge2$. Suppose $\mathcal M$ contains the image of an injective homogeneous map $\bigoplus_{j=1}^mT(-d_j)\to T^m$, where $d_j\ge0$ and $\max_j(d_j-1)\le D$. If $\sum_jd_j>0$, let $\mathcal F'=\bigoplus_{j=1}^m\bigoplus_{b=0}^{d_j-1}T(-b)$, and let $\beta_1(D)$ be the first constant of a $D$-exact decomposition of $\mathcal F'/w_1\mathcal F'$. Then $\reg\mathcal M\le\beta_1(D)$; otherwise, if $\sum_jd_j=0$, then $\mathcal M=T^m$.
\end{lemma}
\begin{proof}
Let $\mathcal L$ be the image of the injective map. If all $d_j$ vanish, its matrix is constant and invertible. Otherwise, the Hilbert function matches that of the quotient by $w_1$, namely
\[
\operatorname{HS}_{T^m/\mathcal L}(t) = \sum_j\frac{1-t^{d_j}}{(1-t)^r} = \operatorname{HS}_{\mathcal F'/w_1\mathcal F'}(t).
\]

Put $N=T^m/\mathcal M$.  If $N=0$, then $\reg\mathcal M=0$ and the bound holds. Assume $N\ne0$.  By combining the generators of $\mathcal L$ with generators of $\mathcal M$ in degrees at most $D$ (and adding a redundant degree-$D$ element if needed), Liang's SPLIT construction and EXACT refinement \cite[Theorem 20 and Lemma 24]{Liang} yield a $0$-exact decomposition $\mathcal{Q}$ of the normal-form space of $N$. Since $N$ is generated in degree zero, we have $b_0(N)=1+\deg\mathcal Q$, which also holds with constant relations since the Hilbert function and $q=0$ fix the constants. The comparison in \cite[Lemma 28]{Liang} provides a $D$-exact decomposition $\mathcal P$ sharing the Hilbert function of $T^m/\mathcal L$ such that $\deg\mathcal Q\le\deg\mathcal P$, where its first constant is $\beta_1(D)$ by \cite[Lemma 26]{Liang}.

Denoting by $l'$ the largest basis weight of $\mathcal F'$, its quotient by $w_1$ matches its Hilbert polynomial in degrees $\ge l'-r+2$.   Since $\beta_1(D)\ge D+\operatorname{rank}\mathcal F'>D\ge l'$, the Hilbert function of the comparison space equals its Hilbert polynomial at each $j\ge\beta_1(D)$. Every positive-dimensional cone has start degree $<\beta_1(D)$, so at each such $j$ its Hilbert function also equals its Hilbert polynomial. The sum of the positive-dimensional cone polynomials is the Hilbert polynomial of the comparison space. Any zero-dimensional cone of degree at least $\beta_1(D)$ would add a positive contribution at its degree, contradicting this equality. Thus $1+\deg\mathcal P=\beta_1(D)$ and $b_0(N)\le\beta_1(D)$. Finally, Nagel's bound \cite[Theorem 4.9]{Nagel} gives $\reg N\le b_0(N)-1$, and the desired regularity follows from the short exact sequence $0\to\mathcal M\to T^m\to N\to0$.
\end{proof}

\begin{lemma}[Bounded lattice weights]
\label{lem:boundedweights-macaulay}
Let $T=K[w_1,\ldots,w_r]$ with $r\ge3$, and let $d_j\ge0$ have a positive sum $s$. For any integer $\sigma\ge\max\{3,s,\max_j d_j\}$, define the module
\[
\mathcal{F}'=\bigoplus_j\bigoplus_{b=0}^{d_j-1}T(-b).
\]
 Then the first constant $\beta_1(\sigma)$ of its $\sigma$-exact decomposition satisfies the upper bound
\[
\beta_1(\sigma)\le2\sigma^{2^{r-2}}.
\]
\end{lemma}
\begin{proof}
Put $t=r-1$ and $D=\sigma$. As a graded free module over $K[w_2,\ldots,w_r]$, the module $\mathcal F'/w_1\mathcal F'$ has $s$ basis elements.
Their weights are $0,\ldots,d_j-1$, for each $j$. Let $W$ be the sum of these weights.
We apply the weighted decomposition construction in \cite[Lemma~34]{Liang}.

At a stage with $k$ variables remaining, let $h_1,\ldots,h_m$ be the residual basis, with weights $\omega_i=\deg h_i$.
We construct the $k$-dimensional cones by assigning them consecutive initial degrees
\[
 b_{k+1},b_{k+1}+1,\ldots,b_{k+1}+m-1,
\]
and set
\[
 b_k=b_{k+1}+m.
\]
Choose a remaining variable $w$. For each $i$, split off the cone with vertex
\[
 h_iw^{b_{k+1}+i-1-\omega_i}.
\]
The lower powers of $w$ give the residual basis over the other $k-1$ variables, with weights
\[
 \omega_i,\omega_i+1,\ldots,b_{k+1}+i-2.
\]
Initially all weights are at most $D$. After this step, every residual weight is at most $b_{k+1}+m-2=b_k-2$.
Hence, at every stage, we have $\omega_i\le b_{k+1}$. In particular, each exponent $b_{k+1}+i-1-\omega_i$ is nonnegative.
 Thus the construction produces a $D$-exact decomposition with consecutive initial degrees for its positive-dimensional cones.

Starting with $b_{t+1}=D$, the first two terms are
\[
 b_t=D+s,\qquad b_{t-1}=(s+1)D+\frac{s(s+1)}2-W.
\]
Since $s\le D$ and $W\ge0$, we get
\[
 \frac{b_t^2}{2}-b_{t-1}=\frac{D^2}{2}-D-\frac s2+W\ge\frac{D(D-3)}2\ge0.
\]
 At any later stage, bounding the number of residual basis elements by $mb_{k+1} + m(m-1)/2$ gives
 \[
 b_{k-1}\le\frac{b_k^2-b_{k+1}^2+b_k+b_{k+1}}2.
\]
Iterating this inequality using $b_k\le b_{k+1}^2/2$ when $b_{k+1}\ge2$, we arrive at
\[
\beta_1(\sigma) \le 2\left(\frac{D+s}{2}\right)^{2^{r-2}} \le 2\sigma^{2^{r-2}}.
\]
For $r=3$, only the initial half-square estimate is needed, completing the proof.
\end{proof}

\paragraph{\textbf{Frobenius presentations.}}
We recover the normal algebra from a free lattice and a common denominator via the Frobenius criterion given below.

\begin{lemma}[Frobenius criterion]
\label{lem:index-Frobenius}
Let $P=K[z_1,\ldots,z_r]$ with $\deg z_i=1$, and let $C$ be a normal graded domain finite over $P$. Suppose that
\[
U=\bigoplus_{j=1}^mPu_j\subseteq C,\qquad \deg u_j=d_j\ge0,
\]
is a free submodule of full rank. Let $b\in P$ be a homogeneous element of degree $a\ge1$ such that $bC\subseteq U$. If $\operatorname{char}K=p>0$ and $q$ is a power of $p$ with $q>ma$, then
\[
C=\{v\in b^{-1}U:v^q\in b^{-1}U\}=b^{-1}\{h\in U:h^q\in b^{q-1}U\}.
\]
\end{lemma}
\begin{proof}
We first prove that every element satisfying the Frobenius condition belongs to $C$. Let $\mathfrak Q$ be an arbitrary height-one prime of $C$.
Since $C$ is finite and integral over $P$, the prime $\mathfrak Q$ lies above a height-one prime $\mathfrak p$ of $P$.
Moreover, the ramification index of $\mathfrak Q$ over $\mathfrak p$ is at most $[\operatorname{Frac}C:\operatorname{Frac}P]=m$.
This bound also holds in the inseparable case. Consequently,
\[
 0\le v_{\mathfrak Q}(b)
 \le m\operatorname{ord}_{\mathfrak p}(b)
 \le ma<q.
\]
Now take an element $v\in b^{-1}U$ such that $v^q\in b^{-1}U$. Since $bC\subseteq U\subseteq C$ and $v^q\in b^{-1}U$, we have $bv^q\in U\subseteq C$.
Therefore, for every height-one prime $\mathfrak Q$ of $C$,
\[
 v_{\mathfrak Q}(b)+qv_{\mathfrak Q}(v)
 =v_{\mathfrak Q}(bv^q)\ge0.
\]
It follows that $qv_{\mathfrak Q}(v) \ge -v_{\mathfrak Q}(b)>-q$.
Since $v_{\mathfrak Q}(v)$ is an integer, we must have $v_{\mathfrak Q}(v)\ge0$.
Thus $v$ belongs to every height-one localization of $C$. The normality of $C$ yields
\[
 v\in\bigcap_{\operatorname{ht}\mathfrak Q=1}C_{\mathfrak Q} =C.
\]
Conversely, let $v\in C$. The assumption $bC\subseteq U$ gives $v\in b^{-1}U$. Since $C$ is a ring, we also have $v^q\in C$, and hence
$v^q\in b^{-1}U$. Therefore,
\[
C=\{v\in b^{-1}U:v^q\in b^{-1}U\}.
\]

Finally, writing $v$ as  $v=h/b$ with $h\in U$, we see that $v^q\in b^{-1}U$ is equivalent to $h^q/b^q\in b^{-1}U$, or equivalently, $h^q\in b^{q-1}U$. Hence
\[
 C=b^{-1}\{h\in U:h^q\in b^{q-1}U\}.
\]
This is the valuation argument of \cite[Theorem~1.1]{SS}; it does not require $U$ to be a subalgebra.
\end{proof}

The next lemma replaces the coordinate equations by a square system while preserving the degrees of its columns.
\begin{lemma}[Compression by rank]\label{lem:rank-congruence-compression}
Let $T=K[w_1,\ldots,w_r]$, where $r\ge2$, let $D\in T$ be nonzero and homogeneous, and let $B\in T^{N\times e}$. Suppose that every nonzero row $\lambda$ has an integer weight $t_\lambda\ge0$ and that $\deg B_{\lambda j}=a_j-t_\lambda$ whenever the entry is nonzero. Set $E=\{v\in T^e:Bv\in DT^N\}$. After a purely transcendental
extension of $K$, there is an $e\times e$ matrix $B'$ such that
\[
 E=\{v\in T^e:B'v\in DT^e\},\qquad
 \deg B'_{ij}=a_j\quad(B'_{ij}\ne0).
\]
The equality is understood after the same field extension.
\end{lemma}
\begin{proof}
We can assume that $K$ is infinite. To make all entries in each column have the same degree, we choose a linear form $\ell$ which divides none of the irreducible factors of $D$, and replace row $\lambda$ by $\ell^{t_\lambda}B_\lambda$. Denote the resulting
matrix by $\widetilde B$. Its $j$th column is homogeneous of degree $a_j$. For every irreducible factor $\phi$ of $D$, the element
$\ell$ is a unit in the discrete valuation ring (DVR) $O=T_{(\phi)}$, since $\ell$ is coprime to $D$. Thus the $O$-row module
\[
 M_\phi=O^e+\sum_\lambda O\,(B_\lambda/D)
\]
remains unchanged when $B$ is replaced by $\widetilde B$. Since $O$ is a DVR, $M_\phi$ is free of rank $e$, and its finite-length quotient $M_\phi/O^e$ can be generated by at most $e$ elements.

Next, we reduce the number of rows from $N$ to $e$ by taking generic linear combinations of the rows of $\widetilde B$. More precisely, we adjoin independent indeterminates $u_{i\lambda}$ to $K$ and let row $i$ of $B'$ be $\sum_\lambda u_{i\lambda}\widetilde B_\lambda$.
For each irreducible factor $\phi$ of $D$, the images of the rows of $B'/D$ span the vector space
\[
 (M_\phi/O^e)/\phi(M_\phi/O^e),
\]
since a corresponding maximal minor is a nonzero polynomial in the $u_{i\lambda}$. Nakayama's lemma therefore shows that these
classes generate $M_\phi/O^e$. The same holds simultaneously for the finitely many factors of $D$. Moreover, a purely transcendental extension preserves the irreducibility of these factors, so the generic coefficients can be chosen to satisfy all the required conditions simultaneously.

Thus the systems agree at every height-one prime dividing $D$, while at every other height-one prime the divisibility conditions are automatic. Since $T$ is a UFD, divisibility by $D$ can be tested at the height-one primes corresponding to its irreducible factors.
Therefore, the two systems define the same submodule of $T^e$. Finally, taking constant linear combinations of rows preserves the column degrees. This completes the proof.
\end{proof}

\begin{lemma}[A Frobenius presentation]
\label{lem:frobenius-presentation}
Use the notation of Lemma~\ref{lem:index-Frobenius} in positive characteristic, and choose a power $q$ of the characteristic
with $q>ma$. After a purely transcendental extension $K\subset K'$, put
\[
 T=(K')^q[w_1,\ldots,w_r],\qquad w_i=z_i^q,
 \qquad \deg_Tw_i=1.
\]
There are a nonzero $D\in T_a$ and an $m\times m$ matrix $B'$ of nonzero determinant, whose $j$th column has degree $d_j$, such that the homogeneous map
\[
 \mathcal A=[B'\mid-DI_m]:
 \left(\bigoplus_{j=1}^mT(-d_j)\right)\oplus T(-a)^m
 \longrightarrow T^m
\]
satisfies
\begin{equation}\label{eq:frobenius-presentation-reg}
 \reg_PC=\reg_T\ker\mathcal A-a.
\end{equation}
Its image is generated in degrees at most $\max\{a,\max_jd_j\}$ and contains the image of the injective map defined by $B'$.
\end{lemma}
\begin{proof}
Since $u_j\in C$ and $bC\subseteq U$, we can write $bu_j^q$ as a $P$-linear combination of the basis elements $u_i$ of $U$.
The coordinate in position $i$ is homogeneous of degree $a+qd_j-d_i$. Put $k=K^q$, and let $L\subseteq K$ be generated over $k$ by
the coefficients of these coordinates and of $b$. The extension $L/k$ is finite, including when $K$ is imperfect. Set
\[
 T_0=k[w_1,\ldots,w_r],\qquad
 U_L=\bigoplus_{i=1}^mL[z_1,\ldots,z_r]u_i.
\]
If the $\lambda_\nu$ form a $k$-basis of $L$, a $T_0$-basis of $U_L$ is
\[
 \lambda_\nu z^\alpha u_i,
 \qquad 0\le\alpha_j<q,\quad 1\le i\le m.
\]
Let $D=b(z)^q\in T_0$, and let $B$ be the coordinate matrix of the $m$ columns $bu_j^q$ in this basis. For $h=\sum_jh_j(z)u_j\in U$, put $v_j(w)=h_j(z)^q$.
By the preceding criterion, we have
\[
 h\in bC
 \quad\Longleftrightarrow\quad
 b\sum_jv_j(w)u_j^q\in b^qU_L
 \quad\Longleftrightarrow\quad Bv\in DT_0^N.
\]
For the first equivalence we initially test membership in $b^qU$. All coordinates lie in $L[z_1,\ldots,z_r]$, and $b^qU\cap U_L=b^qU_L$ by faithful flatness of $L\subseteq K$. Thus coefficientwise Frobenius identifies $bC$, over the graded ring isomorphism $P\simeq T_0$, with
\[
 \mathcal E_0=\left\{v\in\bigoplus_jT_0(-d_j):
                         Bv\in DT_0^N\right\}.
\]
A nonzero row belonging to $z^\alpha u_i$ has integer weight
\[
 t=\frac{|\alpha|+d_i-a}{q}\ge0,
 \qquad \deg B_{\lambda j}=d_j-t
 \quad(B_{\lambda j}\ne0).
\]
The integrality of $t$ follows from the homogeneity of the corresponding coordinate. Moreover, since $q>ma$ and $|\alpha|+d_i\ge0$, we have $t\ge0$.
The columns are independent over $\operatorname{Frac}T_0$. Indeed, the coefficients of a relation have $q$th roots in $\operatorname{Frac}P$. After cancelling $b$, the relation is the $q$th power of a relation among the basis elements $u_j$. Since the $u_j$ are linearly independent over $\operatorname{Frac}P$, all coefficients must be zero.

Apply Lemma~\ref{lem:rank-congruence-compression}. Use independent indeterminates $\mathbf U$ over $k$ for the coefficients of its
auxiliary linear form and row combinations. We may also impose the nonvanishing of a suitable minor to ensure that the resulting matrix has full column rank.
Over $T=T_0\otimes_k k(\mathbf U)$, the resulting square matrix $B'$ has the required column degrees, nonzero determinant, and congruence kernel
\[
 \mathcal E=\left\{v\in\bigoplus_jT(-d_j):B'v\in DT^m\right\}.
\]
Projection onto the first coordinates identifies $\ker\mathcal A$ with $\mathcal E$.

Take independent indeterminates $\mathbf V$ over $K$, set $K'=K(\mathbf V)$, and identify $\mathbf U=\mathbf V^q$. Then $(K')^q=K^q(\mathbf V^q)=k(\mathbf U)$.
Frobenius is a graded ring isomorphism $K'[z_1,\ldots,z_r]\simeq T$ when the $w_i$ have degree one.
Under this isomorphism, $\mathcal E$ corresponds to $b(C\otimes_KK')$. In particular, an element of degree $l$
in $\mathcal E$ corresponds to an element $h\in b(C\otimes_KK')$ of degree $l$, so that $h/b$ has degree $l-a$. Therefore
\[
 \reg_T\mathcal E
 =\reg_{K'[z_1,\ldots,z_r]}(C\otimes_KK')+a
 =\reg_PC+a.
\]
Finally, the matrix $\mathcal A$ has $m$ rows and $2m$ columns, with column degrees $d_1,\ldots,d_m$ and
$a,\ldots,a$. Hence, its size and column degrees are independent of the choice of $q$.
\end{proof}

\begin{lemma}[Transfer to characteristic zero]
\label{lem:char-zero-transfer}
Fix integers $d\ge2$, $a\ge1$, and $c_1,\ldots,c_m\ge0$. Suppose a bound $B(d,\mathbf c,a)$ holds in positive characteristic for every geometrically integral normal graded algebra $V$, finite over $P=K[z_1,\ldots,z_d]$, that contains a full-rank free  lattice
\[
 U=\bigoplus_{j=1}^mPu_j,\qquad \deg u_j=c_j,
\]
and has a nonzero $\Delta\in P_a$ with $\Delta V\subseteq U$. Then the same bound holds in characteristic zero.
For this transfer, it suffices to know the positive-characteristic bound over finite fields.
\end{lemma}
\begin{proof}
Let the data be defined over a characteristic-zero field $K$. Choose finite presentations of $V$ as an algebra and as a $P$-module, together with the lattice inclusion, the element $\Delta$, and a minimal graded $P$-free resolution of $V$. The inclusion $\Delta V\subseteq U$ can be expressed
by finitely many identities involving the $P$-module generators of $V$. Spread all these data over a finitely generated integral $\mathbb Z$-subalgebra $T\subset K$.

After localizing $T$, we may assume that the model $V_T$ and the quotient $V_T/U_T$ are $T$-flat, that the lattice map is injective, and that the spread complex is a free resolution of $V_T$. Indeed, generic freeness gives the required flatness, while a further localization makes the finitely generated kernel and homology modules vanish. We also retain a nonzero coefficient of $\Delta$ as a unit in $T$. Thus, the identities expressing $\Delta V\subseteq U$ remain valid after passing to any fiber.

The generic fiber is geometrically integral and normal:  geometric integrality descends along field extensions, and normality descends from the base change to $K$ and implies geometric normality in characteristic zero.
After shrinking $\operatorname{Spec}T$ further, we may assume that all fibers are geometrically integral and normal,
by \cite[9.7.7(iv) and 12.1.6(iv)]{EGA3}. This follows because the complement of the geometrically normal locus has constructible image missing the generic point, whose closure can be removed.

Flatness of $V_T/U_T$ preserves the injectivity of the lattice map in each fiber. Moreover, $\Delta_tV_t\subseteq U_t$ ensures that
the lattice has full rank in each fiber. Since a nonzero coefficient of $\Delta$ remains a unit, its image $\Delta_t$ is nonzero and still has degree $a$.
The free resolution remains exact after specialization, since its free terms, $V_T$, and hence its successive syzygies
are $T$-flat. Its matrices retain positive polynomial degree, so the fiber resolutions are minimal and have the same graded shifts.
Consequently, the regularity $\reg_PV$ is unchanged. We can therefore apply the positive-characteristic
bound to any closed fiber of $T$. Since $T$ is a finitely generated $\mathbb Z$-algebra, we may choose such a fiber with finite
residue field. This proves the desired bound in characteristic zero.
\end{proof}

\paragraph{\textbf{The normalization bound.}}
Bounding a common denominator by the lattice degrees gives a regularity bound.
\begin{lemma}[Regularity from a bounded lattice]
\label{lem:normal-fiber-lattice-regularity}
Let $V$ be a geometrically integral normal graded algebra in
nonnegative degrees, with $V_0=K$, finite over
$P=K[z_1,\ldots,z_r]$, where $r\ge3$.
Suppose $V$ contains a full-rank free lattice
\[
 U=\bigoplus_jP(-d_j),\qquad d_j\ge0,
 \qquad\sum_jd_j\le\sigma
\]
for an integer $\sigma\ge3$.
Then
\[
 \reg_PV\le2\sigma^{2^{r-2}}+1-\sigma.
\]
\end{lemma}
\begin{proof}
If $V=U$, then its regularity is at most $\sigma$, which gives the bound. Otherwise, put $Q=V/U\ne0$. Since $V$ is normal and finite over $P$, it is
the integral closure of $P$ in $\Frac V$. Therefore, $V$ is the intersection of its localizations at height-one primes of $P$, so it is reflexive as a $P$-module.
The depth lemma shows that $Q$ has no associated primes of height at least two. For each height-one prime $(\phi)$ of $P$, let
\[
 \lambda_\phi=\operatorname{length}_{P_{(\phi)}}Q_{(\phi)}.
\]
The elementary divisor theorem gives $\phi^{\lambda_\phi}V_{(\phi)}\subseteq U_{(\phi)}$. Taking the product over all height-one primes, we obtain
\[
 b_0V\subseteq U,\qquad
 b_0=\prod_\phi\phi^{\lambda_\phi}.
\]
 Indeed, the image $b_0Q$ vanishes at every height-one localization. If it were nonzero, an associated prime of $b_0Q\subseteq Q$ would have height one by the depth argument above, a contradiction. The factors are homogeneous, and multiplicity associativity gives
\[
 \deg b_0=e(Q)=\sum_jd_j-e_1(V)\le\sigma
\]
by Lemma~\ref{lem:normal-finite-e1}. Multiplying $b_0$ by a suitable power of $z_1$, we obtain a homogeneous element $b\in P_\sigma$ satisfying
$bV\subseteq U$.

In positive characteristic, apply Lemma~\ref{lem:frobenius-presentation} with $C=V$ and $a=\sigma$. The image of the resulting matrix is generated in degrees at most $\sigma$ and contains a full-rank free lattice with basis degrees $d_j$. Lemmas~\ref{lem:free-lattice-exact} and
\ref{lem:boundedweights-macaulay} bound its regularity by $2\sigma^{2^{r-2}}$.
The kernel sequence and \eqref{eq:frobenius-presentation-reg} give
\[
 \reg_PV\le2\sigma^{2^{r-2}}+1-\sigma.
\]
Lemma~\ref{lem:char-zero-transfer}, applied to the lattice $U$ and the denominator $b$, extends this bound to characteristic zero.
\end{proof}

\begin{proof}[Proof of Theorem~\ref{thm:normalization}]
After a purely transcendental extension of the base field, we may assume that $A$ admits a linear Noether normalization $P\subseteq A$.
Regularity is preserved, and normalization commutes with this extension:
\[
 \overline{A\otimes_KK(\mathbf U)}=C\otimes_KK(\mathbf U).
\]
The rank of $C$ over $P$ is $e$. If $e=1$, the normality of $P$ implies that $C=A=P$.

Suppose $e\ge2$. By Lemma~\ref{lem:bounded-field-lattice}, applied with $F=A$, $C$ contains a full-rank free lattice with $d_j\le j$ and $\sum_jd_j\le\binom e2$.
Furthermore, since $C$ is contained in $\operatorname{Frac}A$ and $A$ is geometrically integral, so is $C$. In particular,
 $C_0=K$.
Apply Lemma~\ref{lem:normal-fiber-lattice-regularity} with
\[
 p=2^{r-2},\qquad \sigma=\max\left\{3,\binom e2\right\}.
\]
For $e\ge3$, we have $\sigma=\binom e2\le e^2/2$. Thus,
\[
 \reg C\le2\sigma^p+1-\sigma
 \le2\left(\frac{e^2}{2}\right)^p
 \le e^{2p}.
\]
For $e=2$, we have $\sigma = 3$, so the bound follows from $\reg C\le2\cdot3^p-2\le  4^p = e^{2p}$, since $p\ge2$.
\end{proof}

\subsection{Normalized Rees algebras}
\label{subsec:rees-regularity}
The source variables allow us to obtain a sharper lattice bound. A monic presentation then realizes the lattice quotient as a kernel over a polynomial ring with one fewer variable.

\paragraph{\textbf{The source lattice.}}

Let $I$ be generated in degree $\rho\ge2$  and assume $n\ge3$. Give $A=R[It]$ the total grading $\deg x_i=1$ and $\deg t=1-\rho$. Its normalization is denoted by $C$.
Put
\begin{equation}\label{eq:source-constants}
 m=n-1,\qquad E=\sum_{j=0}^{n-1}\rho^j,\qquad
 L=m(\rho-1)+1,\qquad
 D=\frac{m\rho E}{2}-\binom m2.
\end{equation}
Note that $m\rho E$ is always even. Indeed, if $m\rho$ is odd, then both $m$ and $\rho$ are odd. Hence $n=m+1$ is even, and consequently $E$ is even.

\begin{lemma}\label{lem:source-lattice}
After a purely transcendental extension of $K$, there is a linear Noether normalization $P\subseteq A$ containing $x_n$ such that $\operatorname{Frac}C/\operatorname{Frac}P$ is separable. Writing $e=\operatorname{rank}_P C$, we have $e\le E$. Moreover, $C$ contains a full-rank free lattice
\[
 U=\bigoplus_{i=0}^{e-1}Pu_i,\qquad u_0=1,\qquad d_i=\deg u_i,
\]
whose basis degrees, in increasing order, satisfy
\begin{equation}\label{eq:source-lattice-degrees}
 0\le d_i\le\min\{i,L\},\qquad \delta:=\sum_i d_i\le D.
\end{equation}
If $\ell=\max_i d_i$, every integer from $0$ to $\ell$ occurs among the $d_i$. There is also a nonzero homogeneous $b_0\in P$
such that $b_0C\subseteq U$ and $\deg b_0\le\delta$.
\end{lemma}

\begin{proof}
The multiplicity estimate $e(A)\le E$ follows from \cite[Theorem~3.2.4]{Lynn}. Indeed,
\[
 \dim_K A_j \le\sum_{q=0}^j \binom{j+(\rho-1)q+n-1}{n-1}.
\]
The expression on the right has leading coefficient $E/n!$. Comparing Hilbert polynomials gives the desired estimate. Since a linear Noether normalization has rank $e(A)$ and normalization preserves multiplicity, we obtain $e\le E$.

Since $x_n$ is a nonzero linear form in the domain $A$, generic linear forms can be chosen to extend it to a homogeneous system of parameters. Thus, we may choose
\[
 P=K[x_n,z_0,\ldots,z_m]\subseteq A.
\]
We may also require that the matrix of coefficients of $x_1,\ldots,x_m$ in $z_0,\ldots,z_m$ has rank $m$. After constant row operations, subtraction of multiples of $x_n$, and a change of source coordinates, the parameters can be written as
\begin{equation}\label{eq:source-parameters}
 z_0=t f_0(x),\qquad z_i=x_i+t f_i(x)\quad(1\le i\le m),
\end{equation}
where the $f_i\in I_\rho$ and $f_0\ne0$. These conditions define nonempty open subsets of the space of coefficients. The degree-one elements of $A$ generate the rational function field $K(x_1,\ldots,x_n,t)$, so their differentials span its module of K\"ahler differentials over $K$. Since $dx_n\ne0$, generic choices of the remaining parameters extend $dx_n$ to a basis of this module. Thus, the corresponding fraction-field extension is separable. A purely transcendental extension of $K$ is enough to ensure that all these conditions hold, even when $K$ is finite.

Set $F=\operatorname{Frac}P$ and $H=\operatorname{Frac}C$. In the generic fiber, $x_n$ belongs to $F$, and
\begin{equation}\label{eq:source-affine-CI}
 H\simeq F[X_1,\ldots,X_m,T]/ (Tf_0-z_0,\ X_i+Tf_i-z_i\ (1\le i\le m)).
\end{equation}
Indeed, $Tf_0=z_0\ne0$ implies that $f_0$ is invertible, and $R[t,1/f_0]=A[1/f_0]$. Hence, the quotient is the field $A\otimes_PF$, with no additional affine components, and is a zero-dimensional complete intersection in $m+1$ variables. Here each $f_i=f_i(X,x_n)$ has degree at most $\rho$ in $X$.

We now use the perfect Bezoutian tensor for an affine complete intersection \cite[Section~4, Theorem~4.7]{JFMR}. This duality does not require separability. The cited result applies after extending $F$ to an algebraic closure, and perfectness descends by faithful flatness. Let $\mathsf A$ be the $(m+1)\times m$ constant matrix whose zeroth row is zero and whose remaining rows form the identity matrix, and write $f=(f_0,\ldots,f_m)^{\mathsf t}$. Taking $T$ as the last variable in the divided differences, the tensor is represented by
\[
 \Delta=\det[\mathsf A_j+T\Delta_jf(X,Y)\mid f(Y)]
       \quad\text{in }H\otimes_FH.
\]
It is independent of the right-hand variable $S$. Multiplication by the right-hand unit $S$ replaces the last column by $z-\mathsf A Y$. The resulting representative has right-hand source degree at most $m(\rho-1)+1=L$.  Contracting this tensor against $H^\vee$ in the left factor maps onto the right factor by perfectness. Its image is contained in the $F$-span of the right-hand source monomials of degree at most $L$. Hence those monomials span $H$.

Choose a monomial basis adapted to the filtration by source degree, starting with $1$, and denote its elements by $u_i$. These elements are homogeneous in $R\subseteq A$, and their $P$-span is a graded free lattice. If two consecutive filtration spaces agree before reaching $H$, the smaller space is stable under multiplication by every source variable. Since it contains $1$, it must equal $H$, a contradiction. Hence, there are no gaps among the degrees, and $d_i\le i$. Moreover, the spanning estimate gives $d_i\le L$.

To estimate the sum of the basis degrees, we  multiply $S\Delta$ by the left-hand unit $T^{-m}$. Since $T^{-1}=f_0(X)/z_0$, a representative of the resulting perfect tensor is
\[
 \det\bigl[
  \Delta_j f(X,Y)+(f_0(X)/z_0)\mathsf A_j
       \mid z-\mathsf A Y\bigr].
\]
Its total degree in $X,Y$ is at most $m\rho$. Indeed, a term involving $k<m$ columns from the second summand has degree at most $k\rho+(m-k)(\rho-1)+1\le m\rho$. When $k=m$, the apparent term of degree $m\rho+1$ vanishes because $\det[\mathsf A\mid\mathsf A Y]=0$.
Since the basis is adapted to the source-degree filtration, every source monomial of degree at most $v$ is an $F$-linear combination of the $u_i$ with $d_i\le v$. Hence reduction with respect to this basis does not increase the source degree. Thus, the tensor matrix has zero $(i,j)$-entry whenever $d_i+d_j>m\rho$. Since its determinant is nonzero, there is a permutation $\pi$ such that
\[
 d_i+d_{\pi(i)}\le m\rho\quad\text{for every }i.
\]
Put $c=m-1$. Since $d_i\le L=m\rho-c$, we have
\[
 m\rho-d_i-d_{\pi(i)}
 \ge(c-d_i)_++(c-d_{\pi(i)})_+.
\]
Summing over all $i$ gives
\[
 \delta\le\frac{em\rho}{2}-\sum_{i=0}^{e-1}(c-i)_+.
\]
The right-hand side is increasing in $e$, since its successive difference is $m\rho/2-(c-e)_+\ge1$. Using $e\le E$ and $E\ge c$, we obtain $\delta\le D$.

Finally, applying the index-divisor construction in the proof of Lemma~\ref{lem:normal-fiber-lattice-regularity}, we obtain $b_0C\subseteq U$ with
$\deg b_0=\delta-e_1(C)\le\delta$. If $C=U$, we simply take $b_0=1$.
\end{proof}

\paragraph{\textbf{A monic presentation.}}
We now use the source lattice from Lemma~\ref{lem:source-lattice}.
Write
\[
\begin{gathered}
P=B[z],\qquad \dim B=n,\qquad
U=\bigoplus_{i=0}^{e-1}Pu_i\subseteq C,\\
d_i=\deg u_i,\qquad \delta=\sum_{i=0}^{e-1} d_i
\end{gathered}
\]
We choose $z$ appropriately below, while retaining the bounds $\delta\le D$ and $d_i\le L$.

\begin{lemma}\label{lem:monic-defect}
After a purely transcendental field extension and a linear change of coordinates in $P$, the module $Q=C/U$ is finite and torsion-free over
$B$. If $Q\ne0$, there exists a homogeneous polynomial $b\in B[z]$, monic in $z$, such that
\[
 a:=\deg b=\operatorname{rank}_BQ=\delta-e_1(C),\qquad bC\subseteq U.
\]
In positive characteristic  and when $Q\ne0$, the matrix in Lemma~\ref{lem:frobenius-presentation} can be chosen so that, after the coefficientwise Frobenius identification, it defines a map
\[
 G:\bigoplus_iP(-d_i)\longrightarrow P^e,
 \qquad
 \mathcal E_G=\{v:Gv\in bP^e\}\simeq bC,
\]
and $\det G$ is monic in $z$.
\end{lemma}
\begin{proof}
The extension $\Frac C/\Frac P$ is separable by Lemma~\ref{lem:source-lattice}. Thus $C$ is finite \'etale over a nonempty homogeneous open subset of $\operatorname{Spec}P$.
Let $b_0\in P$ be the denominator supplied by that lemma. Choose a line through the origin whose generic point lies in this open subset and does not lie in the zero set of $b_0$. Take this line as the $z$-axis. After scaling, we may assume that $b_0$ is monic in $z$. Therefore, $Q$ is finite over $B$.

Since $C$ and $U$ are reflexive $P$-modules, the depth lemma shows that $\operatorname{depth}Q_{\mathfrak p}\ge1$ whenever $\operatorname{ht}\mathfrak p\ge2$. On the other hand, $Q$ is a torsion $P$-module, so every associated prime of $Q$ has positive height. Therefore, every associated prime of $Q$ has height one. Each such prime contains $b_0$ and contracts to zero in $B$, since $b_0$ is monic in $z$. Thus $Q$ is torsion-free over $B$. Comparing Hilbert series gives
\[
 \operatorname{rank}_B Q=e_P(Q)=\delta-e_1(C)\le \delta.
\]
If this rank is zero, then $Q=0$. Otherwise, let $b(T)$ be the characteristic polynomial of multiplication by $z$ on $Q\otimes_B\Frac B$. Its coefficients lie in $\Frac B$ and are integral over the integrally closed ring $B$, so they belong to $B$. Grading shows that the coefficient of $T^{a-j}$ has degree $j$. Thus $b=b(z)$ is homogeneous of degree $a$ and monic in $z$. By the Cayley--Hamilton theorem, $b$ annihilates $Q\otimes_B\Frac B$. Since $Q$ is torsion-free over $B$, we obtain $bQ=0$. Consequently, $bC\subseteq U$.

 For the remaining assertion, suppose the characteristic is positive and $Q\ne0$. Choose a power $q$ of the characteristic such that $q>ea$.  At the generic point of the $z$-axis, we have $U=C$, and $C/P$ is \'etale. By the fact that relative Frobenius is an isomorphism for \'etale algebras \cite[Tag~0EBS]{Stacks}, the elements $u_i^q$ form a basis there. Since $b$ is a unit at this point, the coordinate matrix of the elements $bu_i^q$ with respect to the basis $u_i$ is invertible there.

We next verify that the Frobenius expansion in Lemma~\ref{lem:frobenius-presentation} preserves this rank condition.
Write $B=K[y_1,\ldots,y_n]$. Set $y_1^q=\cdots=y_n^q=0$ in the expanded matrix and retain the rows whose residue exponents in all the $y_j$ are zero. These rows recover the preceding coordinate matrix after setting $y_1=\cdots=y_n=0$. Hence, a relation among their columns over $K^q(z^q)$ would give a relation among the original columns over $K(z)$. Therefore, the retained matrix has a nonzero $e\times e$ minor at the generic point of the $z$-axis.

Choose the row-homogenizing linear form nonzero on that axis and coprime to $b$. Its powers are units at the generic point of the axis, so row homogenization preserves this rank. Generic constant row combinations can preserve the resulting minor and the finitely many DVR conditions in Lemma~\ref{lem:rank-congruence-compression}. These polynomial nonvanishing conditions hold simultaneously after a purely transcendental extension. The determinant of the compressed matrix has degree $\delta$ and restricts to a nonzero multiple of $z^\delta$. Scaling one row makes it monic in $z$.
The coefficient extension and coefficient roots are as in Lemma~\ref{lem:frobenius-presentation}, so imperfect fields are allowed.
\end{proof}

\begin{lemma}\label{lem:tor-krylov}
Work in positive characteristic with the notation of Lemma~\ref{lem:monic-defect}. Transport all modules through the coefficientwise Frobenius Work in positive characteristic with $Q\ne0$ and with the notation of Lemma~\ref{lem:monic-defect}. Transport all modules through the coefficientwise Frobenius ring isomorphism of Lemma~\ref{lem:frobenius-presentation}, and retain the notation $P=B[z]$, $b$, $U$, and $Q=C/U$ for the transported objects. Put $M=\operatorname{coker}G$. Then $M$ is free over $B$, with a homogeneous basis
\[
 z^jv_i,\qquad 0\le j<d_i,\qquad \deg v_i=0,
\]
and there is an isomorphism of graded $B$-modules
\[
 Q\simeq(0:_M b).
\]
Put $k=\delta-a=e_1(C)$. If $k>0$, the image $bM$ contains a full-rank free $B$-submodule with a basis
\[
 bz^jv_i,\qquad 0\le j<\lambda_i,\qquad
 0\le\lambda_i\le d_i,\qquad \sum_i\lambda_i=k.
\]
In particular, the sum of the degrees of these basis elements is at most $k(a+(L-1)/2)$.
\end{lemma}
\begin{proof}
In column $i$ of $G$, the coefficient of $z^{d_i}$ is constant. The matrix of these leading coefficients is invertible because $\det G$ is monic in $z$. A constant row change makes the leading terms $z^{d_i}e_i$, while all other terms in column $i$ have smaller $z$-degree. Thus, with respect to a monomial order
that compares $z$-degree first, the columns form a monic module Gr\"obner basis over $B$. Indeed, their leading terms have distinct positions, so no further reductions are needed. Division gives the stated $B$-basis of $M$.

Put $F_0=\bigoplus_iP(-d_i)$. Tensoring the exact sequence
\[
 0\longrightarrow F_0\xrightarrow{G}P^e
 \longrightarrow M\longrightarrow0
\]
with $P/bP$ gives
\[
 \operatorname{Tor}_1^P(M,P/bP)\simeq\mathcal E_G/bF_0.
\]
Indeed, $bF_0\subseteq\mathcal E_G$, and the kernel of the induced map $F_0/bF_0\to(P/bP)^e$ consists precisely of the classes of those $v$ satisfying $Gv\in bP^e$.

Under the coefficientwise Frobenius identification, $\mathcal E_G$ corresponds to $bC\subseteq U$, while $bF_0$ corresponds to $bU$. Indeed, for $h=bu$ with $u\in U$, taking the $q$th powers of the coordinates corresponds to multiplying the transported vector of $u$ by the transported polynomial $b$. Therefore,
\[
 \mathcal E_G/bF_0\simeq bC/bU\simeq Q(-a).
\]
On the other hand, the resolution
$0\to P(-a)\xrightarrow{b}P\to P/bP\to0$ gives
\[
 \operatorname{Tor}_1^P(M,P/bP)\simeq(0:_M b)(-a).
\]
Comparing these two descriptions of the same Tor module and cancelling the common degree shift yields
\[
 Q\simeq(0:_M b).
\]
Since $M$ has $B$-rank $\delta$ and $Q$ has $B$-rank $a$, it follows
that
\[
 \operatorname{rank}_B(bM)=\delta-a=k.
\]

Over $F=\Frac B$, consider the map $F[z]^e\to bM\otimes_BF$ sending $e_i$ to $bv_i$. Its kernel contains the columns of $G$. Its initial submodule is
$\bigoplus_i z^{\lambda_i}F[z]e_i$ for some integers $0\le\lambda_i\le d_i$. The standard monomials therefore give an $F$-basis of $bM\otimes_BF$.
Consequently, their images form a $B$-linearly independent set and generate a full-rank free $B$-submodule of $bM$. The dimension of $bM\otimes_BF$ gives
$\sum_i\lambda_i=k$. The sum of the degrees of these basis elements is
\[
 \sum_i\sum_{j=0}^{\lambda_i-1}(a+j)=ak+\sum_i\binom{\lambda_i}{2}\le k\left(a+\frac{L-1}{2}\right).\qedhere
\]
\end{proof}

\paragraph{\textbf{The multiplication image.}}
We now incorporate the prefix condition into the cone comparison. For a finite multiset $\Omega$ of nonnegative weights, start the construction at an integer $A\ge\max\Omega$. If its sorted weights are $\omega_1,\ldots,\omega_h$, the next endpoint is $A+h$, and the residual multiset consists of the intervals
\[
 [\omega_i,A+i-2]\qquad(1\le i\le h),
\]
with multiplicities. Empty intervals are omitted. Iterating this procedure gives the exact cone constants of the weighted free module.
If $M_j=\sum_{\omega\in\Omega}\binom{\omega}{j}$, one step is described by
\begin{equation}\label{eq:cone-moment-recurrence}
 A'=A+M_0,\qquad
 M'_j=\binom{A'}{j+2}-\binom{A}{j+2}-M_{j+1}.
\end{equation}
This formula follows by summing the binomial coefficients over the residual intervals.

We next introduce an order on finite multisets. Write their elements in increasing order, including multiplicities, as $\omega_1,\ldots,\omega_h$ and
$\widetilde\omega_1,\ldots,\widetilde\omega_{\widetilde h}$. We say that the second multiset dominates the first if $\widetilde h\ge h$ and
$\widetilde\omega_i\le\omega_i$ for $1\le i\le h$.
Equivalently, it contains at least as many weights at or below every threshold. If, in addition, $\widetilde A\ge A$, one cone step preserves this order: for each $i\le h$,
\[
 [\omega_i,A+i-2]
 \subseteq
 [\widetilde\omega_i,\widetilde A+i-2],
\]
while the additional weights contribute further intervals. Consequently, $\widetilde A+\widetilde h\ge A+h$, so the next starting
degree is no smaller. Repeating the same argument shows that this inequality persists at every subsequent cone step.
\begin{lemma}\label{lem:packed-cones}
	Let $1\le a<\delta$ and $k=\delta-a$, let $0\le\lambda_i\le L$ satisfy $\sum_i\lambda_i=k$, where $L\ge3$. Set $A=a+L-1$. For each degree $c=a+j$ with
	$0\le j<\lambda_i$, let the comparison string be $0,\ldots,c-1$. Every cone endpoint is bounded above by the corresponding endpoint
	obtained by packing $k$ into strings of length $L$ and at most one shorter string. If $S$ is the number of weights in this packed
	comparison and the polynomial ring has $n\ge3$ variables, then its first cone constant is at most
	\[
	2\left(\frac{A+S}{2}\right)^{2^{n-2}},\qquad S\le k\left(a+\frac{L-1}{2}\right).
	\]
This majorization also holds for any larger integer $k$.
\end{lemma}
\begin{proof}
For every integer $v\ge0$, the number of comparison weights at most $v$
is
\[
\sum_i F_v(\lambda_i),\qquad
F_v(\lambda)=\sum_{j=0}^{\lambda-1}\min\{a+j,v+1\}.
\]
Since the successive differences of $F_v$ are nondecreasing, $F_v$ is discretely convex. Therefore, if $0<\lambda_i\le\lambda_j<L$, replacing $(\lambda_i,\lambda_j)$ by $(\lambda_i-1,\lambda_j+1)$ does not
decrease this sum.  Iterating this operation yields a partition of $k$ into strings of length $L$ and at most one shorter string.
Thus, the packed comparison has at least as many weights as the original comparison below every threshold.   By the multiset order defined above, threshold domination alone gives the cone-step comparison. It follows that every cone endpoint is bounded above by the corresponding endpoint for the packed comparison.

 For a larger integer $k'$, add singleton prefixes to the original collection to reach $k'$ column degrees. This only increases the threshold counts. Packing this enlarged collection again weakly increases every count.  Hence the majorization remains valid when $k$ is replaced by any larger integer.

	Let $W$ be the sum of the packed comparison weights. The first two endpoints are
	\[
	T=A+S,\qquad H=(S+1)A+\frac{S(S+1)}{2}-W.
	\]
	Hence
	\[
	\frac{T^2}{2}-H=\frac{A^2}{2}-A-\frac{S}{2}+W.
	\]
	If $a\ge2$, each degree $c\ge a$ contributes $\binom c2\ge(a-1)c/2$. Therefore, $W\ge(a-1)S/2$, and  we obtain
	\[
	\frac{T^2}{2}-H
	\ge \frac{A^2}{2}-A+\frac{a-2}{2}S
	\ge0.
	\]
	If $a=1$, each full prefix of length $L$ contributes
	\[
	\sum_{c=1}^L\left(\binom c2-\frac c2\right)
	=\frac{L(L+1)(2L-5)}{12}\ge0
	\]
	to $W-S/2$. The remaining shorter prefix, if present, contributes at least $-1/2$. Thus, the difference $\frac{T^2}{2}-H$ is at least
	\[
	\frac{L^2}{2}-L-\frac12>0.
	\]
In either case, we have $H\le T^2/2$. The remaining half-square estimates from Lemma~\ref{lem:boundedweights-macaulay} therefore yield the desired bound.
	
Finally, let $\widetilde\lambda_i$ be the packed lengths. Since $\widetilde\lambda_i\le L$ and $\sum_i\widetilde\lambda_i=k$, we have
	\[
	\sum_i\binom{\widetilde\lambda_i}{2}
	\le\frac{L-1}{2}\sum_i\widetilde\lambda_i
	=\frac{(L-1)k}{2}.
	\]
	Consequently,
	\[
	S=\sum_i\sum_{j=0}^{\widetilde\lambda_i-1}(a+j)
	=ak+\sum_i\binom{\widetilde\lambda_i}{2}
	\le k\left(a+\frac{L-1}{2}\right),
	\]
	which proves the asserted estimate.
\end{proof}

\begin{proposition}\label{prop:tor-bound}
	For the normalized Rees algebra $C$ and the source-lattice bounds $D$ and $L$,
let $p=2^{n-2}$ and
	\[
	J=\frac{(D+(L+1)/2)^2}{8}+\frac{L-1}{4}.
	\]
	Then, over every field,
	\begin{equation}\label{eq:tor-regularity-bound}
		\reg C\le2J^p.
	\end{equation}
	More generally, if $e_1(C)\le K$, then
	\begin{equation}\label{eq:tor-conditional-bound}
		\reg C\le2\left(
		\frac{D+L-1+K(D+(L-1)/2)}2\right)^p.
	\end{equation}
\end{proposition}
\begin{proof}
 First suppose that the characteristic is positive. If $Q=0$, then $C=U$ and $\reg C\le L$. If $Q\ne0$ and $k=e_1(C)=0$, Lemma~\ref{lem:tor-krylov} gives $\operatorname{rank}_B(bM)=0$. Since $bM$ is a submodule of the free $B$-module $M$, we have $bM=0$ and $Q\simeq M$. Thus $\reg_B Q\le L-1$ and $\reg C\le L$. In either case, $D\ge m\rho^n/2\ge\max\{8,L\}$ ensures that the two displayed upper bounds are at least $L$.
Therefore, we may assume that $a,k\ge1$.

Consider the map $M(-a)\xrightarrow{b}M$. In the displayed $B$-basis of $M$, a row of weight $\omega$ can be
multiplied by a fixed nonzero linear form raised to the power $\omega$. This embeds the target in $B^\delta$ without changing the image as an
	abstract graded module. A generic constant projection from $B^\delta$ to $B^k$ is injective on the image over
	$\Frac B$, and hence is injective over $B$.
	Let $N\subseteq B^k$ be its image.
	It is generated in degrees at most $A=a+L-1$ and contains the
	prefix lattice of Lemma~\ref{lem:tor-krylov}.
	
	Apply Lemma~\ref{lem:free-lattice-exact} to $N$, followed by Lemma~\ref{lem:packed-cones}. With $h=(L-1)/2$, the packed comparison
	gives
	\[
	\begin{aligned} A+S&\le a+2h+(\delta-a)(a+h)=a(\delta+1-a-h)+h\delta+2h\\
		&\le\frac{(\delta+1+h)^2}{4}+h\le2J.
	\end{aligned}
	\]
	Consequently,  Lemma~\ref{lem:packed-cones} yields  $\reg_B N\le2J^p$. The kernel sequence gives
	\[
	\reg_B Q\le\max\{L-1,\reg_B N+1-a\}
	\le\max\{L-1,\reg_BN\}.
	\]
	Since $Q$ is finite over $B$, the irrelevant ideals of $P$ and $B$
	have the same radical on $Q$. Thus,
	$\reg_PQ=\reg_BQ$.
	The exact sequence $0\to U\to C\to Q\to0$ now yields
	\eqref{eq:tor-regularity-bound}.
	
	If $k=e_1(C)\le K$, we apply the same packed estimate, using the
	bounds $a\le D$ and $k\le K$. Then
	\[
	A+S\le D+L-1+K\left(D+\frac{L-1}{2}\right),
	\]
	which gives \eqref{eq:tor-conditional-bound}.

Finally, consider the case where the base field has characteristic zero. We use the spreading argument of
Lemma~\ref{lem:char-zero-transfer}.  Spread the lattice $U$, the inclusion $b_0C\subseteq U$, and the multiplication table in the basis $u_i$ after inverting $b_0$. Over the generic field the trace pairing has nonzero determinant by separability. Clear its denominators and invert in the coefficient base a nonzero coefficient of the resulting numerator, as well as a nonzero coefficient of $b_0$. In the fibers the $u_i$ still form a basis over the fraction field of the specialized polynomial ring, and the specialized trace determinant remains nonzero. Thus the fibers are generically separable, and their lattice shifts and bounds are those of the original lattice. The fixed resolution preserves the Hilbert series, and hence preserves both $e_1$ and regularity.
We may therefore choose a positive-characteristic fiber to which the bounds already established apply. Transferring
these bounds back to the original characteristic-zero algebra completes the proof.
\end{proof}

\subsection{The first Hilbert coefficient and explicit bounds}
\label{subsec:hilbert-explicit}

We retain the notation of \eqref{eq:source-constants} and put $p=2^{n-2}$.  An upper  bound  on the first Hilbert coefficient  then yields the smaller exponent in characteristic zero.

\begin{lemma}[Genus under specialization]\label{lem:genus-specialization}
	Let $V$ be an excellent discrete valuation ring with algebraically  closed residue field $K$ of characteristic zero. Let
	$\mathcal X\to\operatorname{Spec}V$ be a proper flat curve whose  generic fiber is smooth and geometrically connected of genus $g$.
	If $Z_1,\ldots,Z_q$ are distinct irreducible components of its reduced special fiber, then
	\[
	\sum_{i=1}^q g(Z_i^\nu)\le g,
	\]
	where $Z_i^\nu$ denotes the normalization of $Z_i$.
\end{lemma}
\begin{proof}
	The generic fiber is geometrically integral, so flatness implies that $\mathcal X$ is integral. After a finite extension of the fraction field and the corresponding  extension of discrete valuation rings, the semistable reduction theorem gives  a semistable model of the generic fiber
	\cite[Tag~0CDN]{Stacks}. By resolving its singularities, we may further assume that this model is regular. Indeed, resolving local equations of the
		form $xy=\pi^a$ introduces only chains of rational curves  \cite[Tag~0CDC]{Stacks}. The residue field remains $K$, and the reduced special fiber of $\mathcal X$ is unchanged.

Take the closure of the graph of the generic-fiber isomorphism between the regular model obtained from the semistable model and the base change of $\mathcal X$. A resolution of this integral excellent surface gives a common regular model $\mathcal W$ by \cite[Tag~0BGP]{Stacks}. Its proper birational morphism to the regular model is a sequence of point blowups by \cite[Tag~0C5R]{Stacks}. Every new special-fiber component in this process and in the resolution of the semistable model is rational. The strict transforms preserve the function fields of the old components. The connected nodal fiber of the semistable model has arithmetic genus $g$ by flatness. Its normalization sequence expresses $g$ as the sum of the genera of its normalizations and the first Betti number of its dual graph. Hence $\sum_{F\subseteq\mathcal W_k}g(F^\nu)\le g$.

	For each $i$,  we choose an irreducible component $F_i$ of $\mathcal W_k$	surjecting onto $Z_i$.
Since the components $Z_1,\ldots,Z_q$ are distinct, the components $F_1,\ldots,F_q$ are also distinct.	The induced morphisms
$F_i^\nu\longrightarrow Z_i^\nu$	are nonconstant  morphisms of smooth projective curves, hence finite.	Applying the Riemann--Hurwitz formula in characteristic zero yields $g(Z_i^\nu)\le g(F_i^\nu)$. Summing these inequalities over all $i$, we obtain
	\[
	\sum_{i=1}^q g(Z_i^\nu)
	\le \sum_{i=1}^q g(F_i^\nu)
	\le \sum_{F\subseteq\mathcal W_k}g(F^\nu)
	\le g.
	\]
	This completes the proof.
\end{proof}

\begin{lemma}\label{lem:charzero-e1}
	Suppose $\operatorname{char}K=0$. Then
	\[
	e_1(C)\le \binom n2(\rho-1)\rho^{n-2}+1+\rho+\cdots+\rho^{n-2}.
	\]
\end{lemma}
\begin{proof}
We may replace $K$ by its algebraic closure, since this preserves normality and the Hilbert polynomial. Set $X=\operatorname{Proj}C$.
A linear Noether	normalization gives a finite map $\pi:X\to\mathbf P^n$ with $\mathcal O_X(1)=\pi^*\mathcal O_{\mathbf P^n}(1)$. Thus
$\mathcal O_X(1)$ is invertible, even if $C$ is not standard graded.
	
We first show that the preimage $D_\Lambda=\pi^{-1}(\Lambda)$ of a general line $\Lambda\subseteq\mathbf P^n$	is a smooth,
possibly disconnected curve. Let $G$ be the Grassmannian of lines in $\mathbf P^n$,  and we consider the incidence scheme of pairs $(x,\Lambda)\in X\times G$ satisfying $\pi(x)\in\Lambda$. The projection from this incidence scheme to $X$ is a $\mathbf P^{n-1}$-bundle. Since $X$ is normal, $\dim \operatorname{Sing}X \le n-2$, so the part of the incidence scheme lying over	$\operatorname{Sing}X$ has dimension at most $2n-3<\dim G$.
	The generic fiber over $G$ therefore lies completely inside the smooth locus $X_{\mathrm{reg}}$, where the	incidence scheme is regular.
	Therefore,  the generic fiber over $G$ is regular and hence smooth in characteristic zero \cite[Tag~038S]{Stacks}.
By properness, $D_\Lambda$ is also smooth for general $\Lambda$.
	
	The $n-1$ linear equations defining $\Lambda$ form a regular sequence near $D_\Lambda$, since  $X$ is regular there and $D_\Lambda$ has codimension $n-1$. Elsewhere on $X$, at least one equation is a unit, so the Koszul complex is exact everywhere.
Computing Euler characteristics from this complex gives
\begin{equation}\label{eq:section-genus}
	\begin{aligned}
		\chi(\mathcal O_{D_\Lambda}(v))&=e(C)(v+1)-e_1(C),\\
		e_1(C)&=e(C)+p_a(D_\Lambda)-1\textcolor{red}{.}
	\end{aligned}
\end{equation}
Note that no global Cohen--Macaulay assumption is needed.

Choose $0\ne f_0\in I_\rho$. Every component of $D_\Lambda$ meets the open set $x_nf_0\ne0$. Indeed, the image of its complement under $\pi$ has dimension at most $n-1$;  a general line  $\Lambda$ is not contained in this image, while every component of $D_\Lambda$ surjects onto $\Lambda$.
On this open set,  the normalization of $X$ agrees with the source chart
	\[
	K[X_1,\ldots,X_{n-1},T,1/f_0(X)],\qquad
	X_i=x_i/x_n,\quad T=t x_n^{\rho-1}.
	\]
In these coordinates, the equations of $D_\Lambda$ become
	\[
	h_j(X,T)=\ell_j(X)+T f_j(X)=0\qquad(1\le j<n),
	\]
	where $\ell_j$ is affine linear and $\deg f_j\le\rho$.
	Thus distinct components of $D_\Lambda$ are birational to distinct generically smooth components of this system.
	
Next, homogenize the $h_j$ on  $Y=\mathbf P^{n-1}\times\mathbf P^1$  to obtain sections $H_j$ of $\mathcal O_Y(\rho,1)$.
Choose general sections $G_j$ of the same line bundle and consider  the family defined by $H_j+\varepsilon G_j=0$ over $K[[\varepsilon]]$.
	Its generic fiberis a smooth,   geometrically connected complete intersection curve, whose genus is given by the adjunction formula as
	\[
	g=1+\frac{(n+1)(n-2)\rho^{n-1}-n(n-1)\rho^{n-2}}2.
	\]
	
		Let $\mathcal X$ be the schematic closure of this generic curve in $Y_{K[[\varepsilon]]}$, which	is proper and flat over $K[[\varepsilon]]$.
	At the generic point of each source component identified above, the   Jacobian of the $h_j$ has rank $n-1$.  By the Jacobian criterion, the perturbed intersection is	smooth over $K[[\varepsilon]]$ near that point,  and agrees with $\mathcal X$ there,  since flatness excludes vertical components. Thus each source component occurs as a distinct reduced component of the special fiber of $\mathcal X$.
	
Finally,  we decompose $D_\Lambda$ into its connected components as $D_\Lambda=D_1\sqcup\cdots\sqcup D_q$. Lemma~\ref{lem:genus-specialization} gives $\sum_{i=1}^q g(D_i)\le g$,  and  $p_a(D_\Lambda)=\sum_{i=1}^q g(D_i)-q+1\le g$.
Substituting $e(C)\le E$ into \eqref{eq:section-genus}, we obtain
	\[
	e_1(C)\le E+g-1=\binom n2(\rho-1)\rho^{n-2}+1+\rho+\cdots+\rho^{n-2}.
	\]
	This proves the assertion.
\end{proof}

Set
\[
H_0=\frac{D+L-1+K_0\bigl(D+(L-1)/2\bigr)}{2},
\]
where $K_0 = \binom{n}{2}(\rho - 1)\rho^{n-2} + 1 + \rho + \cdots + \rho^{n-2}$.
Proposition 3.14 and Lemma 3.16 give
\begin{equation}\label{eq:reg-C-bound}
	\operatorname{reg} C\le 2H_0^p.
\end{equation}
Moreover, the transfer construction in Proposition 3.14	preserves $e_1(C)$, so the same upper bound $K_0$ remains valid.

\paragraph{\textbf{Explicit estimates.}}
\label{subsec:explicit-bound}
 To express  the bounds in terms of $n$ and $\rho$,  we first compute the exact cone recurrence needed at $(n, \rho) = (4, 2)$.

 \begin{lemma}\label{lem:critical-cones}
	Suppose $D=42$ and $L=4$. For each integer $1\le a\le41$, we write $42-a=4q+r$, where $0\le r<4$, and form the comparison strings as in Lemma~\ref{lem:packed-cones}, using
	$q$ blocks of length $4$ and one remainder block. Starting at $A=a+3$, the third cone endpoint is at most
	\[
	K_{*}=6\,874\,078\,799<288^4.
	\]
\end{lemma}
\begin{proof}
Write $42-a=4q+r$, where $0\le r<4$. The column degrees are $a,a+1,a+2,a+3$, each repeated $q$ times, followed by $a,\ldots,a+r-1$.
Let $S$, $W$, and $V$ denote, respectively, the total cardinality, the sum of the weights, and the sum of the second binomial moments of these comparison strings.
	Then
	\[
	\begin{aligned}
		A&=a+3,\qquad S=q(4a+6)+ra+\binom r2,\\
		W&=q\left[\binom{a+4}{3}-\binom a3\right]
		+\binom{a+r}{3}-\binom a3,\\
		V&=q\left[\binom{a+4}{4}-\binom a4\right]
		+\binom{a+r}{4}-\binom a4.
	\end{aligned}
	\]
Applying \eqref{eq:cone-moment-recurrence} three times gives
	\begin{equation}\label{eq:critical-endpoint}
		\begin{aligned}
			T&=A+S,\qquad
			H=(S+1)A+\frac{S(S+1)}2-W,\\
			K(a)&=\frac{H^2-T^2+H+T}{2}-\binom T3+\binom A3+V.
		\end{aligned}
	\end{equation}
Thus, $K(a)$ is the third cone endpoint.

We now determine the maximum of $K(a)$ by means of a Bernstein-coefficient certificate, avoiding the need to compute and list all $41$ values individually.
For $0\le i<4$, let $r_i$ be the least nonnegative residue of $42-i$ modulo $4$.
Define $K_i(x)$ by substituting $r=r_i$ and $q=(42-x-r_i)/4$ into \eqref{eq:critical-endpoint}, and set
	\[
	P_i(x)=K_{i+1}(x+1)-K_i(x),
	\]
	with subscripts read modulo $4$. Thus, $P_i(a)=K(a+1)-K(a)$ whenever $a\equiv i\pmod4$. Each $P_i$ is a polynomial of degree at most $7$.

For each $i$, we determine intervals on which $P_i$ is positive or negative by expressing the corresponding polynomial in the Bernstein basis of degree $7$.
 More precisely, the table below lists intervals $I_+$ and $I_-$ on which all the degree-$7$ Bernstein coefficients of $P_i$ and $-P_i$, respectively, are positive.
The last two columns give the corresponding smallest Bernstein coefficients.
	\[
	\begin{array}{c|c|c|r|r}
		i&I_+&I_-&b_+&b_-\\\hline
		0&[4,20]&[24,40]&27\,096\,840&11\,776\,590\\
		1&[1,17]&[21,37]&38\,607\,376&30\,961\,424\\
		2&[2,18]&[22,38]&72\,920\,025&77\,972\,442\\
		3&[3,19]&[23,39]&84\,578\,440&35\,864\,340
	\end{array}
	\]
The entries are obtained directly from \eqref{eq:critical-endpoint}. Indeed, if $P_i(\alpha+(\beta-\alpha)u)=\sum_{j=0}^7c_ju^j$,
	then its Bernstein coefficient in position $\ell$ is
	\[
	\sum_{j=0}^{\ell}c_j\frac{\binom{\ell}{j}}{\binom{7}{j}}, \text{\ where\ } 0\le\ell\le7.
	\]
Since the Bernstein basis polynomials are nonnegative on $[0,1]$ and sum to $1$, positivity of all Bernstein coefficients implies strict positivity of the
	corresponding polynomial throughout the interval. Consequently, the intervals $I_+$ show that $K(a+1)-K(a)>0$ for every integer $1\le a\le20$,
	while the intervals $I_-$ show that $K(a+1)-K(a)<0$ for every integer $21\le a\le40$.
	Therefore, $K(a)$ increases up to $a=21$ and decreases thereafter. Hence,
\[
	K(a)\le K(21)=6\,874\,078\,799<6\,879\,707\,136=288^4. \qedhere
	\]
\end{proof}

\begin{lemma}\label{lem:explicit-regularity}
Let $n\ge4$ and $\rho\ge2$. With the notation above,
\[
 \reg C\le M_n(\rho)
 =\left(\frac{(n-1)^2\rho^{2n}}8\right)^{2^{n-2}}.
\]
If $\operatorname{char}K=0$, then we also have
\[
 \reg C\le
 \left(\frac{n(n-1)^2\rho^{2n-1}}4\right)^{2^{n-2}}.
\]
\end{lemma}
\begin{proof} 	Put $X=m\rho^n$. We first consider the case $\rho\ge3$. By the estimates $D\le3X/4$, $(L+1)/2\le X/8$, and $L-1\le X/16$,
	we obtain
	\[
	8J\le\frac{49X^2}{64}+\frac X8<\frac{4X^2}{5},
	\]
	where the last inequality follows from $X\ge243$. Since $p\ge4$ and $2(4/5)^p<1$, it follows that
	\begin{equation}\label{eq:J-clean}
		2J^p<\left(\frac{X^2}{8}\right)^p.
	\end{equation}
	We next consider $\rho=2$ and $m\ge4$. In this case, $X=4mp$, $L=m+1$, and $D=X-m(m+1)/2$. Writing $v=(m^2-2)/2$, we have
	\[
	8J=(X-v)^2+2m.
	\]
We claim that
\[
 X^2-8J=(m^2-2)X-v^2-2m\ge3mX.
\]
Indeed, the inequalities
\[
 m^2-3m-2\ge\frac m2,\qquad
 X\ge m^3,\qquad
 v^2+2m\le\frac{m^4}{2}
\]
imply
\[
 \begin{aligned}
 X^2-8J
 &= (m^2-3m-2)X
    +3mX-v^2-2m\\
 &\ge \frac{mX}{2}+3mX-\frac{m^4}{2}\ge 3mX.
 \end{aligned}
\]
Therefore,
\[
 \frac{8J}{X^2}\le1-\frac{3}{4p}.
\]
Moreover,
\[
 \left(1-\frac{3}{4p}\right)^{-p}
 \ge1+\frac34+\frac{9(p+1)}{32p}>2.
\]
Thus, we again obtain
\[
 2J^p<\left(\frac{X^2}{8}\right)^p.
\]
This establishes \eqref{eq:J-clean} in both cases. By Proposition~\ref{prop:tor-bound}, we obtain the desired
bound over an arbitrary field, except possibly when
$(n,\rho)=(4,2)$.

In this exceptional case, $D=42$ and $L=4$. If $a,k\ge1$, we enlarge $k$ to $42-a$ in Lemma~\ref{lem:packed-cones} and retain the exact cone
endpoints. By Lemma~\ref{lem:critical-cones}, we have
\[
 \reg C\le K_*<288^4.
\]
If $a=0$ or $k=0$, the desired bound follows directly from $\reg C\le L$. Furthermore, the characteristic-zero transfer
in Proposition~\ref{prop:tor-bound} preserves the same exact cone estimate. Since
\[
 288^4=\left(\frac{(4-1)^2\cdot2^{2\cdot4}}8\right)^{2^{4-2}}=M_4(2),
\]
the first assertion follows.

We now prove the improved bound when $\operatorname{char}K=0$.
Put
\[
 Y=\frac{n(n-1)^2\rho^{2n-1}}4.
\]
We use the following estimates:
\[
 \begin{aligned}
 K_0&\le
 \left[\frac{nm}{2}(\rho-1)+2\right]\rho^{n-2},
 &D&\le\frac{m\rho^{n+1}}{2(\rho-1)},\\
 K_0&\le\frac{nm}{2}\rho^{n-1},
 &L-1&\le m\rho.
 \end{aligned}
\]
Combining these inequalities with $n\ge4$ and $\rho\ge2$, we obtain
\[
 \begin{aligned}
 DK_0&\le
 \left(\frac{nm^2}{4}+m\right)\rho^{2n-1},\\
 \frac{K_0(L-1)}2&\le
 \frac{nm^2}{32}\rho^{2n-1},\\
 D+L-1&\le
 \frac{9m}{64}\rho^{2n-1}.
 \end{aligned}
\]
Therefore,
\[
 \frac{H_0}{Y}
 \le\frac9{16}+\frac{73}{32nm}
 \le\frac{289}{384}.
\]
Since $p\ge4$, we have $2^{1/p}<6/5$. Hence
\[
 2^{1/p}\frac{H_0}{Y}
 \le\frac65\cdot\frac{289}{384}<1,
\]
which gives $2H_0^p<Y^p$.
Finally, equation~\eqref{eq:reg-C-bound} yields
\[
 \reg C\le Y^p
 =\left(\frac{n(n-1)^2\rho^{2n-1}}4\right)^{2^{n-2}},
\]
as desired.
\end{proof}

\subsection{A source restriction}
\begin{lemma}\label{lem:full-defect-section}
	Let $\ell$ be a general source linear form, put $S=R/(\ell)$, and set
	\[
	\mathcal D=C/\ell C,\qquad T=H^0_{\mathfrak m_x}(\mathcal D),\qquad
	E=\operatorname{im}(C\longrightarrow S[t]),\qquad
	C'=\overline{S[(IS)t]}.
	\]
	Then $E=\mathcal D/T$, the module $T$ is concentrated in the zeroth component of the first grading, and $\mathfrak m_xT=0$. Moreover,
	\[
	\reg_xE_{(*,s)}=\reg_R\overline{I^s}-\rho s
	\qquad(s\ge1).
	\]
	The algebra $C'$ is the normalization of $E$. For each $s$, the module $W_{(*,s)}$, where $W=C'/E$, has finite length over $S$.
	If $\ell(I)<n$, then $T=0$.
\end{lemma}
\begin{proof}
The Rees valuations give rise to primary representations of all $\overline{I^s}$ whose centers belong to a fixed finite set. Thus, a general linear form $\ell$ can be chosen to be filter-regular on every $R/\overline{I^s}$ and  to satisfy $IS\ne0$. By Lemma~\ref{lem:normal-saturation},
$H^1_{\mathfrak m_x}(C)_{(u,*)}=0$ for all $u\ge0$. The long exact sequence of local cohomology induced by multiplication by $\ell$, together with the absence of negative first degrees in $C$, shows that $T$ is concentrated in the zeroth component of the first grading. Hence $\mathfrak m_xT=0$.

	Write $J_s=\overline{I^s}$. The kernel of the natural map $J_s/\ell J_s\to S$ is $((J_s:\ell)/J_s)(-1)$.
Since $\ell$ is filter-regular on $R/J_s$, this kernel has finite length over $S$ and is therefore contained in $T$. Conversely, every element of $T$ maps to zero in $S[t]$, since $S[t]$ has no $\mathfrak m_x$-torsion. Hence the kernel is exactly $T$, and $E=\mathcal D/T$.
	
	Each $E_{(*,s)}$ is a nonzero ideal of $S$, shifted by $\rho s$, and contains a nonzero element of shifted degree zero. Its regularity is therefore nonnegative. The exact sequence
	\[
	0\longrightarrow T_s\longrightarrow\mathcal D_s\longrightarrow E_s\longrightarrow0
	\]
	and the fact that $T_s$ is concentrated in shifted degree zero imply
	\[
	\reg_x\mathcal D_s=\reg_xE_s.
	\]
	Indeed, if $T_s\ne0$, then $\reg_xT_s=0\le\reg_xE_s$, and the regularity inequalities for a
	short exact sequence give the asserted equality; the case $T_s=0$ is immediate.
	Since $\ell$ is regular on each ideal strand, this yields
	\[
	\reg_xE_{(*,s)}=\reg_R J_s-\rho s,
	\]
	as desired.
	
	The extension $S[(IS)t]\subseteq E$ is finite and has fraction field $\operatorname{Frac}(S)(t)$, so its normalization is $C'$.
	To establish normality away from the source vertex, write $\ell=\sum a_ix_i$ with algebraically independent coefficients.
	Before passing to $K(\mathbf a)$, eliminating $a_j$ gives
	\[
	(C[\mathbf a]/(\ell))_{x_j}
	\cong C_{x_j}[a_i:i\ne j].
	\]
	The ring on the right is normal, and the elimination identifies its map into the ambient source ring as an inclusion.
Thus, $E_{x_j}$ is normal. Since $E\subseteq C'$ is a finite birational extension, we obtain
	\[
	E_{x_j}=C'_{x_j} \qquad\text{for every }j.
	\]
	Thus $W_{(*,s)}$ is supported only at the source vertex and has finite length over $S$.
	
		Finally, the valuative criterion gives   $(\overline{I^{s+1}}:I)=\overline{I^s}$.
	Multiplication by a set of generators of $I$ therefore embeds $R/\overline{I^s}$ into a direct sum of copies of
	$R/\overline{I^{s+1}}$. Hence the associated prime sets of the normal powers form an ascending chain. If $\ell(I)<n$, McAdam's theorem
	\cite[Theorem~5.4.6]{HS} excludes $\mathfrak m$ from their stable set, and hence from every one of them. Since this ascending chain stabilizes, its union is finite and does not contain $\mathfrak m$. We may therefore choose $\ell$ outside the union, so that it is regular on every $R/\overline{I^s}$.
	 The restriction then has kernel $\ell C$, and consequently $T=0$.
\end{proof}

\subsection{The component of first degree zero}
Suppose $\ell(I)=n$, and write
\[
B=C_{(0,*)},\qquad A_0=K[I_\rho t],\qquad
E_0=E_{(0,*)},\qquad V=C'_{(0,*)}.
\]
The rings $B$ and $V$ are normal graded domains; for example, $B=C\cap\Frac(B)$. Taking first degree zero in a finite Rees extension makes them finite over the corresponding ordinary fiber algebras. We identify $B$ with a subalgebra of $R^{(\rho)}$.
\begin{lemma}\label{lem:zero-face-lattice}
Let $k=\rho^{n-2}$. There exists a linear Noether normalization
	\[
	P_0=K[z_1,\ldots,z_n]\subseteq A_0,\qquad
	r=\operatorname{rank}_{P_0}B\le\rho^{n-1},
	\]
	such that $B$ contains a full-rank free lattice
	$U=\bigoplus_jP_0(-d_j)$ satisfying
	\[
	d_j\ge0,\qquad \sum_jd_j\le k(r-1).
	\]
	For the ideal $J=\ker(B\to E_0)$, there exists a nonzero homogeneous
	$c\in P_0\cap J$ of degree $1\le h\le k$.
Setting  $\theta_i=x_i/\ell$,  we have
	\[
	J=\{v\in B:\theta_iv\in B\text{ for all }i\}.
	\]
\end{lemma}
\begin{proof}
Let $f_1,\ldots,f_n\in I_\rho$  be general elements defining a generically finite map $\varphi:\mathbb P^{n-1}\dashrightarrow\mathbb P^{n-1}$,
and set $z_j=f_jt$. These elements  define  a linear normalization of  $A_0$. Let $F=\Frac(P_0)$ and let $H$ be the fraction field of the
source Veronese algebra. The generic degree $r=[H:F]$ is at most $\rho^{n-1}$ by residual B\'ezout.
	
For a source ratio $\theta_i=x_i/\ell$, let $\Gamma_i$ be the closure of the graph of $(\theta_i,\varphi)$, and let
$Y_i\subseteq\mathbb P^1\times\mathbb P^{n-1}$ be its image. The image is an integral hypersurface with equation
\[
\sum_{j=0}^{a_i}c_{ij}(Y)Z_0^jZ_1^{a_i-j}=0,
\qquad c_{i,a_i}\ne0,\qquad \deg c_{ij}=b_i,
\]
where $a_i=[F(\theta_i):F]$. The generic degree of $\Gamma_i\to Y_i$ is   $\delta_i=r/a_i$. Indeed, the corresponding
	extension of function fields has degree $r/a_i$, since the 	intermediate function field obtained by adjoining $\theta_i$
	has degree $a_i$ over $F$. These function-field degrees agree 	with the generic degrees of the corresponding projective maps:
	adjoining a common nonzero homogeneous parameter of degree one 	recovers the original graded fraction fields.

Let $\alpha$ and $\beta$ be the hyperplane classes from $\mathbb P^1$ and $\mathbb P^{n-1}$, respectively, and use the same
letters for their pullbacks to $\Gamma_i$. Since $[Y_i]=a_i\alpha+b_i\beta$, the projection formula gives
\[
\delta_i b_i=\deg\bigl([\Gamma_i]\alpha\beta^{n-2}\bigr).
\]
Choose a general divisor from the first factor and $n-2$ general hyperplanes from the second. Their pullbacks are basepoint-free
on $\Gamma_i$. Since the complement of the source graph has dimension at most $n-2$, these $n-1$ general divisors
	can be chosen so that their intersection avoids this complement. On the source open, the cuts are one equation of degree one and
$n-2$ equations of degree $\rho$. Residual B\'ezout therefore yields
\[
\delta_i b_i\le\rho^{n-2}=k;
\]
see \cite[Sections~1.4, 2.3 and Theorem~12.3]{Fulton}.
Here the intersection number is understood with multiplicities, so the estimate also accounts for inseparable generic degrees.

	Evaluating at $z_j=f_jt$, define
\[
\zeta_{i,j}=c_{i,a_i}\theta_i^j+
c_{i,a_i-1}\theta_i^{j-1}+\cdots+
c_{i,a_i-j+1}\theta_i
\quad(1\le j<a_i).
\]
These elements have degree $b_i$ and are integral over $P_0$.  To see this, consider an arbitrary valuation of $H$
that is nonnegative on $P_0$. If its value on $\theta_i$ is 	nonnegative, the displayed expression shows that its value on
$\zeta_{i,j}$ is nonnegative. Otherwise, the defining equation 	allows us to rewrite $\zeta_{i,j}$ as a negative sum of terms
involving nonpositive powers of $\theta_i$, which again have 	nonnegative values. Thus, the valuative criterion for integral
dependence shows that $\zeta_{i,j}$ is integral over $P_0$.
This is the classical binary-form basis \cite[Section~2.1]{Wood}. The same argument shows that
$c_{i,a_i}\theta_i$ is integral over $P_0$.   Under the embedding of the source Veronese algebra in $R[t]$ given by $g\mapsto g t^s$ for $\deg g=\rho s$, we have $H\subseteq\Frac (A)$. Since $P_0\subseteq A$, each of these integral elements is integral over $A$ and hence belongs to $C$. Its first degree is zero, so it belongs to $B$. Dividing $c_{i,a_i}\theta_i$ by the nonzero $c_{i,a_i}$ shows that $\Frac (B)$ contains each $\theta_i$. These ratios generate the degree-zero rational function field of the source. Adjoining $z_1$ to that field generates $H$, so $\Frac (B)=H$.

	Set $F_i=F(\theta_1,\ldots,\theta_i)$, $r_i=[F_i:F_{i-1}]$, and $D_i=[F_i:F]$.  At step $i$, the elements $1,\zeta_{i,1},\ldots,\zeta_{i,r_i-1}$
	form a basis over $F_{i-1}$, by their triangular transition matrix from  powers   of $\theta_i$.
Products of these bases yield a homogeneous $F$-basis in $B$. Its total degree is bounded by
\[
	\sum_jd_j\le r\sum_i b_i\frac{r_i-1}{r_i}\le k\sum_i a_i\frac{r_i-1}{r_i}\le k\sum_i(D_i-D_{i-1})=k(r-1),
\]
since $\delta_i=r/a_i$, $\delta_i b_i\le k$,
$a_i\le D_i$, and $D_i=D_{i-1}r_i$.
Taking the $P_0$-span of this homogeneous $F$-basis gives a free $P_0$-submodule of $B$ of rank $r$,
	which is the required full-rank free lattice. This proves the lattice bound.
	
Choosing $\ell$  also  to attain the minimum source value at each Rees valuation of $I$, if $v\in J$,  writing $v=\ell w$ yields
$v_Q(x_iw)\ge v_Q(\ell w)$ for every Rees valuation $v_Q$, so $\theta_iv=x_iw\in B$ by the valuative criterion.
Conversely, if all $x_iv/\ell$ are polynomials, choose an $x_i$ that is coprime to $\ell$. Then $\ell$ divides $v$.
This proves the formula for $J$.
For such an $x_i$,  substituting into  its graph equation  and setting $c=c_{i,a_i}\in P_0$,  multiplying by $\ell^{a_i}$ and
reducing modulo $\ell$ gives $c(f_1,\ldots,f_n)x_i^{a_i}=0$. Thus $\ell$ divides $c(f_1,\ldots,f_n)$, so $c\in J$.
Its degree satisfies $1\le h=b_i\le k$.
\end{proof}

\begin{lemma}\label{lem:fractional-intersection-bound}
	Let $B$ be a normal geometrically integral graded domain finite over
	a standard polynomial ring $P$ in $r\ge3$ variables. Let
	$U=\bigoplus_jPu_j\subseteq B$ be a full-rank free lattice with
	$\deg u_j=d_j\ge0$, and let $b\in P_a$, $a\ge1$, satisfy $bB\subseteq U$.
	For $\theta_0=1,\theta_1,\ldots,\theta_N\in\Frac(B)_0$, set
	\[
	J=\{v\in B:\theta_iv\in B\text{ for every }i\}.
	\]
	If $c\in P_h\cap J$ is nonzero and
	$\Sigma\ge\max\{3,\sum_j(d_j+h),a+h\}$ is an integer, then
	\[
	\reg_PJ\le2\Sigma^{2^{r-2}}+1-\Sigma.
	\]
\end{lemma}
\begin{proof}
	Put $w_j=cu_j$ and choose $\beta\in P_\Sigma$ divisible by $bc$.
	Then $\beta J\subseteq cU\subseteq J$.
	
	In characteristic $p>0$, we  choose a power $q>ma$, where
	$m=\operatorname{rank}_PB$. By
	Lemma~\ref{lem:index-Frobenius}, for every $v\in\Frac(B)$,
	\[
	v\in B\quad\Longleftrightarrow\quad v^q\in b^{-1}U.
	\]
	Indeed, the forward implication follows from $B\subseteq b^{-1}U$, and the reverse implication follows from $bv^q\in U\subseteq B$.
Writing $H=\sum_jg_jw_j\in cU$, we obtain
	\[
	H/\beta\in J\quad\Longleftrightarrow\quad
	b\sum_jg_j^q(\theta_iw_j)^q\in\beta^qU
	\quad\text{for every }i.
	\]
	Since $\theta_iw_j\in B$, these conditions are polynomial congruences in the basis of $U$.
	
	As in Lemma~\ref{lem:frobenius-presentation}, we expand the coordinates over	$K^q[z_1^q,\ldots,z_r^q]$, with
	$\deg z_i^q=1$. The resulting system has nonzero rows indexed by $z^\alpha u_l$, each with weight
	$t=(|\alpha|+d_l-a)/q\ge0$ and $j$th entry of degree $d_j+h-t$. The block $i=0$ has full column rank.
	By Lemma~\ref{lem:rank-congruence-compression},	the system can therefore be reduced to a square
		matrix with nonzero determinant, column degrees $d_j+h$, and common modulus of degree $\Sigma$.
	The resulting kernel presentation represents $\beta J$	under coefficientwise Frobenius. Therefore, Lemmas~\ref{lem:free-lattice-exact} and~\ref{lem:boundedweights-macaulay} yield an upper bound
$2\Sigma^{2^{r-2}}$  for the regularity of the 	corresponding image.
	Combining this estimate with the kernel exact sequence and the degree shift by $-\Sigma$ gives
	\[
	\reg_PJ\le2\Sigma^{2^{r-2}}+1-\Sigma.
	\]
	The coefficient extension and generic coefficient roots are as in Lemma~\ref{lem:frobenius-presentation}; perfectness is unnecessary.
	
	In characteristic zero, write $a_i=c\theta_i\in B$. The module $J$ is the kernel of
	\[
	B\longrightarrow(B/cB)^N,\qquad
	v\longmapsto(a_iv\bmod cB)_i.
	\]
Spread this map, its kernel, the lattice, the	denominators, and a minimal $P$-free resolution of $J$ over a finitely generated integral
$\mathbb Z$-algebra. After shrinking the base, generic freeness ensures that the kernel commutes with specialization and that the free
	resolution remains exact and minimal. As in Lemma~\ref{lem:char-zero-transfer}, we may also assume that
	the fibers of $B$ are geometrically integral and normal, and that $b$ and $c$ remain nonzero.
The positive-characteristic estimate then applies to these fibers with the same degree data and Betti shifts,
		yielding the desired bound.
	\end{proof}
	
\subsection{Transfer to the power strands}

Fix a basis $f_1,\ldots,f_N$ of $I_\rho$ and put $Y=(y_1,\ldots,y_N)$.
	We use the bigraded presentation ring $R[Y]$, with $\deg y_j=(0,1)$ and
	$y_j\mapsto f_jt$ (or its image under restriction). Put
	\[
	Z=2\rho^{2n-3},\qquad \Gamma=2Z^{2^{n-2}}+1-Z.
	\]

\begin{proposition}\label{prop:stabilization-transfer}
	Suppose $\reg C\le\tau$ and $\reg C'\le\eta$, where $\tau,\eta$ are integers with $\tau\ge\eta+1$. If $\ell(I)=n$, suppose also
	$\tau\ge\Gamma$. Then
	\[
	\rstabn(I)\le\max\{\rstabn(IS),\tau\}.
	\]
\end{proposition}
\begin{proof}
 	Let $W=C'/E$. For $u\ge1$, Lemma~\ref{lem:full-defect-section} gives $E_{(u,*)}=(C/\ell C)_{(u,*)}$. Since $\ell$ is $C$-regular and has
	total degree one, $\reg(C/\ell C)=\reg C$. Taking fixed first-degree strands in the bigraded free resolutions of $C/\ell C$ and $C'$ over $R[Y]$ yields
	\[
	\reg_{K[Y]}E_{(u,*)}\le\tau,\qquad
	\reg_{K[Y]}C'_{(u,*)}\le\eta.
	\]
	Indeed, the nonnegative first coordinates of the bigraded Betti shifts ensure that the regularity of each strand is bounded by the total regularity.
The exact sequence
\[
0\longrightarrow E\longrightarrow C'\longrightarrow W\longrightarrow0
\]
 thus yields
	\[
	\reg_{K[Y]}W_{(u,*)}\le\tau-1.
	\]
	If $\ell(I)<n$, the same argument applies to $u=0$, since $T=0$.
In the full-spread case,  Lemmas~\ref{lem:normal-fiber-lattice-regularity}, \ref{lem:zero-face-lattice}, and \ref{lem:fractional-intersection-bound} supply
 the same bound for $u=0$.

	We next show that the source-degree support of $W_s$ is independent of $s\ge\tau$.
 For each $u$, the Hilbert function of $W_{(u,*)}$ equals its Hilbert polynomial for $s\ge\tau$. Since $W_{(u,*)}$ is generated in degrees at most $\tau-1$, a zero component in this range forces all later components to vanish. The polynomial is then identically zero, so the components are either all zero or all nonzero in this range. Since each $W_s$ has finite source-degree support, we obtain
	\[
	\operatorname{end}_x W_s
	\quad\text{is independent of }s\ge\tau.
	\]
		
	The strands $E_s$ and $C'_s$ are nonzero ideals of $S$, with the same shift $\rho s$. Since $\dim S\ge2$, their zeroth source local cohomology modules vanish.  Moreover, $W_s$ has finite length, since its source-degree support is finite. Hence the local-cohomology sequence associated to
\[
0\longrightarrow E_s\longrightarrow C'_s
\longrightarrow W_s\longrightarrow0
\]
gives
\[
\reg_xE_s=
\max\{\reg_xC'_s,\operatorname{end}_xW_s+1\}.
\]
The second term on the right-hand side is therefore constant for all $s\ge\tau$. On the other hand, since $IS$ is nonzero and generated in degree $\rho$,  the regularity of $\overline{(IS)^s}$ has eventual slope $\rho$, and its defect is constant for all $s\ge\rstabn(IS)$.
It follows that the defect of $E_s$ is constant for every
$s\ge\max\{\rstabn(IS),\tau\}$.

	Finally, applying Lemma~\ref{lem:full-defect-section} identifies the corresponding defect of $E_s$ with  that of $\overline{I^s}$. Hence
\[
\rstabn(I)\le\max\{\rstabn(IS),\tau\},
\]
as desired.
\end{proof}

\subsection{Proof of the main results}
\begin{proof}[Proof of Theorem~\ref{thm:main}]
	Proposition~\ref{thm:reduction} reduces the assertion to an ideal generated in degree $\rho=\rho(I)\le d$. If $\rho=1$,
	Proposition~\ref{prop:elementary-cases} applies. Assume $\rho\ge2$. We argue by induction on $n$. Retain the quantities $D,L,J,p$
	of the preceding section, and write
	\[
	F_n(\rho)=2J^{2^{n-2}}.
	\]
	By the definitions of these quantities, $F_n(\rho)$ is increasing in both $n$ and $\rho$. Proposition~\ref{prop:tor-bound} gives $\reg C\le F_n(\rho)$.
	It remains to show that the same bound controls $\rstabn(I)$.
	
	Let $Z=2\rho^{2n-3}$ and $\Gamma=2Z^{2^{n-2}}+1-Z$. We first verify the numerical condition needed in Proposition~\ref{prop:stabilization-transfer}. For $n\ge4$, the inequalities $J>D^2/8$ and $D\ge(n-1)\rho^n/2$ give
	\[
	\frac JZ>\frac{(n-1)^2\rho^3}{64}\ge\frac98>1.
	\]
	For $n=3$, we have $D=\rho^3+\rho^2+\rho-1\ge\rho^2(\rho+1)$, and hence
	\[
	\frac JZ>\frac{\rho(\rho+1)^2}{16}\ge\frac98>1.
	\]
	Thus $J>Z$ for every $n\ge3$, and consequently
	\[
	F_n(\rho)=2J^{2^{n-2}}
	>2Z^{2^{n-2}}
	>\Gamma.
	\]
The restricted normalization has regularity at most $\eta=(1+\rho+\cdots+\rho^{n-2})^{2^{n-1}}\le\Gamma-1$, by Theorem~\ref{thm:normalization} and the multiplicity estimate in Lemma~\ref{lem:source-lattice}. Since $\Gamma$ is an integer, we may take
$\tau=\lfloor F_n(\rho)\rfloor$, which satisfies $\tau\ge\Gamma\ge\eta+1$. Moreover,
$\reg C\le\tau$. We may therefore apply Proposition~\ref{prop:stabilization-transfer}. Proposition~\ref{prop:elementary-cases} gives the two-variable case, and induction together with the monotonicity of $F_n$ gives the required bound
for the stabilization index of the reduced ideal. Taking second-degree strands of a minimal bigraded resolution of $C$ and using Proposition~\ref{prop:tor-bound} also gives $\Bn(I)+1\le F_n(\rho)$.

	For $n\ge4$, except at $(n,\rho)=(4,2)$, inequality~\eqref{eq:J-clean}
	gives $F_n(\rho)\le M_n(\rho)$. At that endpoint, the exact cone calculation instead gives
	\[
	\reg C\le K_*=6\,874\,078\,799<M_4(2).
	\]
	This integer also exceeds $F_3(2)=52\,441/32$ and $\Gamma_4(2)=33\,554\,369$. Thus the induction bound for $n=3$ is also dominated by $K_*$, and the numerical condition $\tau\ge\Gamma_4(2)$ is satisfied.
	The transfer argument with $\tau=K_*$ bounds $\rstabn(I)$ by $K_*$;  the same strand estimate gives the same bound for $\Bn(I)+1$.
	This proves the all-field statement for $n\ge4$, with $\rho$ in place of $d$. For $n=3$,
	\[
	D+\frac{L+1}{2}=\rho^3+\rho^2+2\rho-1<2\rho^3,
	\qquad J<\rho^6.
	\]
	Hence $F_3(\rho)<2\rho^{12}\le\rho^{16}$, which gives the stated bound in this case.
	
		For $n\ge4$ in characteristic zero, let $H_0$ be the quantity  defined before \eqref{eq:reg-C-bound}.  It increases with $n$ and $\rho$.
	For $n\ge4$, $H_0\ge(n-1)\rho^{2n-1}/4>Z$. For $n=3$, $K_0=3\rho^2-2\rho+1>2\rho^2$ and $D\ge\rho^3$, so again
	$H_0>Z$.
	Therefore, using the characteristic-zero estimate for $\reg C$ in place of the general estimate, we obtain $2H_0^p>\Gamma$.
	The same induction with $2H_0^p$ gives
	\[
	\max\{\Bn(I)+1,\rstabn(I)\}\le2H_0^p
	\le\left(\frac{n(n-1)^2\rho^{2n-1}}4\right)^p,
	\]
	where the last inequality is the characteristic-zero numerical estimate. Combining these bounds and applying monotonicity in $\rho\le d$ completes the proof.
	\end{proof}

\begin{proof}[Proof of Corollary~\ref{cor:regularity-powers}]
	Assume first that $n\ge4$ and $d\ge2$. Put $p=2^{n-2}$ and	$t=(n-1)2^n$. Since $n-1\le2^{n-2}\le d^{n-2}$, we have
	\[
	M_n(d)=\left(\frac{(n-1)^2d^{2n}}8\right)^p
	\le8^{-p}d^{4(n-1)p}\le d^t.
	\]
	Theorem~\ref{thm:main} therefore gives $\Bn(I)\le d^t$.	For $n=3$, Theorem~\ref{thm:main} gives $\Bn(I)\le d^{16}$.
Observing that $16=(3-1)2^3 \le (n-1)2^n$ for all $n\ge3$, both cases can be unified as
	\[
\Bn(I)\le d^{(n-1)2^n}.
	\]
	Thus, for $n\ge3$ and $d\ge2$,
	\[
	\reg\overline{I^s}\le\rho(I)s+\Bn(I)
	\le \rho(I)s+d^{(n-1)2^n}\qquad(s\ge1).
	\]
	If $n\le2$ or $d=1$, Proposition~\ref{prop:elementary-cases} gives $\reg\overline{I^s}=\rho(I)s$, so the assertion follows.
\end{proof}

\hspace{-6mm} {\bf Acknowledgment: }
This research is supported by the Natural Science Foundation of Jiangsu Province (No.  BK20221353) and the National Natural Science Foundation of China (No. 12471246).  The second author carried out part of this work while visiting the International Center for Research in Mathematics and Postgraduate Training (ICRTM), affiliated with the Institute of Mathematics at the Vietnam Academy of Science and Technology, through the ICRTM--TWAS Visiting Program. She expresses her sincere gratitude for the financial support and warm hospitality received during this period.

\end{document}